\documentclass[aos,preprint]{imsart}

\RequirePackage[OT1]{fontenc}
\RequirePackage{amsthm,amsmath}
\usepackage{listings}
\usepackage{latexsym}
\DeclareSymbolFont{AMSb}{U}{msb}{m}{n}
\DeclareSymbolFontAlphabet{\mathbb}{AMSb}

\startlocaldefs

\newcommand{\nc}{\newcommand}
\nc{\nt}{\newtheorem}
\nt{defn}{Definition}
\nt{lem}{Lemma}
\nt{pr}{Proposition}
\nt{theorem}{Theorem}
\nt{cor}{Corollary}
\nt{ex}{Example}
\nt{ass}{Assumption}
\nt{step}{Step}
\nt{case}{Case}
\nt{subcase}{Subcase}
\nt{note}{Note}
\nt{remark}{Remark}
\nt{algo}{Algorithm}

\nc{\bd}{\begin{defn}} \nc{\ed}{\end{defn}}
\nc{\blem}{\begin{lem}} \nc{\elem}{\end{lem}}
\nc{\bpr}{\begin{pr}} \nc{\epr}{\end{pr}}
\nc{\bth}{\begin{theorem}} \nc{\eth}{\end{theorem}}
\nc{\bcor}{\begin{cor}} \nc{\ecor}{\end{cor}}
\nc{\bex}{\begin{ex}}  \nc{\eex}{\end{ex}}
\nc{\bass}{\begin{ass}}  \nc{\eass}{\end{ass}}
\nc{\bstep}{\begin{step}}  \nc{\estep}{\end{step}}
\nc{\bcase}{\begin{case}}  \nc{\ecase}{\end{case}}
\nc{\bsubcase}{\begin{subcase}}  \nc{\esubcase}{\end{subcase}}
\nc{\bnote}{\begin{note}}  \nc{\enote}{\end{note}}
\nc{\prf}{{\bf Proof.} }
\nc{\eop}{\hfill $\Box$ \\ \\}
\nc{\argmin}{\mathrm{argmin}}
\nc{\argmax}{\mathrm{argmax}}
\nc{\sgn}{\mathrm{sgn}}
\nc{\Var}{\mathrm{Var}}
\nc{\Cov}{\mathrm{Cov}}
\nc{\bak}{\!\!\!\!\!}
\nc{\IBD}{\mathrm{IBD}}
\nc{\supp}{\mathrm{supp}}
\nc{\dom}{\mathrm{dom}}
\nc{\R}{{\mathbb R}}
\nc{\peq}{\preceq}
\nc{\wt}{\widetilde}
\nc{\Mult}{\mathrm{Mult}}
\nc{\Prob}[1]{\mathbb{P}_{#1}}

\endlocaldefs

\begin{document}
\bibliographystyle{plain}

\begin{frontmatter}
\title{Estimating a probability mass function with unknown labels}
\runtitle{Estimating a probability mass function}

\begin{aug}
\author{\fnms{Dragi Anevski} \thanksref{m1}\ead[label=e1]{dragi@maths.lth.se}},
\author{\fnms{Richard D. Gill} \thanksref{m2}\ead[label=e2]{gill@math.leidenuniv.nl}}
\and
\author{\fnms{Stefan Zohren} \thanksref{m3}
\ead[label=e3]{stefan.zohren@materials.ox.ac.uk}
\ead[label=u1,url]{}}

\runauthor{Anevski, Gill and Zohren}

\affiliation{Lund University\thanksmark{m1}, Leiden University\thanksmark{m2} and University of Oxford\thanksmark{m3} }

\address{Centre for Mathematical Sciences, \\ Lund University, \\
Box 118, 221 00 Lund, \\ Sweden\\
\printead{e1}}

\address{Mathematical Institute, \\ Leiden University\\
Niels Bohrweg 1, 2333 CA Leiden, \\ The Netherlands \\
\printead{e2}}

 \address{
Department of Materials,\\ 
Parks Road, OX1 3PH, Oxford, UK \\
\printead{e3}}
\end{aug}

\begin{keyword}[class=MSC]
\kwd{62G05, 62G20, 65C60, 62P10}
\end{keyword}

%\begin{keyword}
%\kwd{sample}
%\kwd{\LaTeXe}
%\end{keyword}
\runauthor{Anevski, Gill and Zohren}

\begin{center}
{\bf Abstract}
\end{center}
\begin{abstract}
In the context of a species sampling problem we discuss a non-parametric maximum likelihood estimator for the underlying probability mass function. The estimator is known in the computer science literature as the high profile estimator. We prove strong consistency and derive the rates of convergence, for an extended model version of the estimator. We also study a sieved estimator for which similar consistency results are derived. Numerical computation of the sieved estimator is of great  interest for practical problems, such as forensic DNA analysis, and we present a computational algorithm based on the stochastic approximation of the expectation maximisation algorithm. As an interesting byproduct of the numerical analyses we introduce an algorithm for bounded isotonic regression for which we also prove convergence. 
\end{abstract}

\end{frontmatter}

\begin{keyword}[class=MSC]
%\kwd[Primary ]{}
\kwd{62G05, 62G20, 65C60, 62P10}
%\kwd[; secondary ]{}
%\end{keyword}
%\begin{keyword}
%\kwd{}
%\kwd{}
%\end{keyword}
\end{frontmatter}

% AOS,AOAS: If there are supplements please fill:
%\begin{supplement}[id=suppA]
%  \sname{Supplement A}
%  \stitle{Title}
%  \slink[doi]{10.1214/00-AOASXXXXSUPP}
%  \sdatatype{.pdf}" 
%  \sdescription{Some text}
%\end{supplement}

\section{Introduction}

Assume we have a random sample that is drawn from an infinite population of species. The goal of this paper is to, based on the random sample, estimate the unknown relative frequencies of all the species in the population. 

Probably the most well known estimator in the context of species sampling is the \emph{naive estimator}, which is the vector of relative frequencies of the species observed in the sample. The problem of this estimator is that it assigns zero probability to any new species which have not yet been observed in the sample. However, when the relative frequencies are very small it is very likely that when sampling a new element this will be a new, so far unobserved species. Such a situation arises for example in forensic DNA analysis when the Y-STR profile of the suspect is not present in the database. This makes it necessary to go beyond the naive estimator and consider estimators for the unknown relative frequencies of all the species in the population.

The first to have studied problems in this setting  
is apparently Fisher et al.\ \cite{fisher:corbet:williams:1943}, who assumed that the members of each separate species are caught according to separate Poisson processes with different intensities and allowing for the processes to be dependent.

The first to use a non-parametric approach is Good \cite{good:1953}, who presented an approximate formula for the expectation of the population frequency. Good attributes the formula to Alan Turing. His approximation becomes better for larger sample sizes but it is not clear from the results in his paper if the formula is asymptotically correct. As a consequence he is also able to give an estimate of the coverage, the sum of the population frequencies of the species observed in the sample, leading to what is known as the \emph{Good estimator} or Good Turing estimator for the probability mass of the unobserved species,
%\[
%\hat{Q}_n=\frac{1}{n} \sum_i1_{\{N_i=1\}},
%\]
%i.e. 
which is given by
the number of species observed exactly once in a sample devided by the sample size.  Next Good and Toulmin \cite{good:toulmin:1956} study a similar setting but for the case when there is a second sample drawn from the population, which can then be thought of as an enlargement of the original sample. As an application Efron and Thisted \cite{efron:thisted:1976} used the result by Fisher et al.\ \cite{fisher:corbet:williams:1943} and Good and Toulmin \cite{good:toulmin:1956} to estimate the number of words known by Shakespeare based on the observed word frequencies in his works. Later work has been concerned with the bias, confidence intervals as well as asymptotic normality of the Good estimator (e.g.\ \cite{esty,esty2,zhang:zhang:2009}), see also Mao and Lindsay \cite{mao} for an application to DNA analysis in this context.

One sees that the naive estimator and the Good estimator are complementary in the sense that the former gives an estimate for the probability distribution of the already observed species, while the latter gives an estimate for the total probability mass of all unobserved species. One would like to combine both these estimators and extend the tail of the naive estimator over the region of unobserved species. A proposal for such an estimator has been made in \cite{orlitsky:sajama:santhanam:viswanathan:2004a,orlitsky:sajama:santhanam:viswanathan:2004b,orlitsky:sajama:santhanam:viswanathan:2005,acharya:orlitsky:pan:2009} for a similar problem in a computer science setting. In \cite{orlitsky:sajama:santhanam:viswanathan:2004a} they introduced what they call the \emph{high profile estimator} and what we refer to as the \emph{pattern maximum likelihood estimator} (PML) which is explained in detail below. For small models this estimator can be obtained analytically \cite{orlitsky:sajama:santhanam:viswanathan:2004a,acharya:orlitsky:pan:2009} and for bigger models a Monte Carlo expectation maximisation (EM) algorithm was proposed in \cite{orlitsky:sajama:santhanam:viswanathan:2004b}.  In \cite{orlitsky:sajama:santhanam:viswanathan:2005} they have also claimed, without complete proof, consistency results for the PML, and discussed the general problem of modelling and estimation of the distribution over ``large alphabets'' when there is a small sample. Their work has been the main motivation for the research presented here. In particular our goals have been to give a full consistency proof, as well as an extension of their model together with its numerical implementation. % \\

We can state the basic estimation problem of the high profile estimator or PML in a simplified manner as follows: Given $N_1,\ldots,N_K$ a set of absolute frequencies, $N_i$ denoting the number of times a species $i$ is observed, and ordered (by us) in decreasing order. There is another order, provided by Nature, which orders the species in how frequent they are in Nature, modelled by a set of decreasing probabilities  $\theta_1,\theta_2,\ldots$ that sum to one, where $\theta_{\alpha}$ denotes how frequent the $\alpha$th most frequent species  is. We can view our data $N$ as an ordering of an underlying  data set $X_{\alpha_1},\ldots,X_{\alpha_k}$ (for some indices $\alpha_i, i=1,\ldots,n$). There is an unobserved map, which takes the order provided by us to the order provided by Nature, which we denote by $\chi$ and which is a bijection. We will derive the likelihood for $\theta$ based on the data $N$ for this problem, and define the PML of $\theta$ as the maximizer of that likelihood under the assumptions $\theta_1\geq \theta_2\geq \ldots, \sum \theta_i=1$. However, typically, and with high probability, the PML $\hat{\theta}$ will not exist in the above model.

Therefore, besides the above described, \emph{basic model}, we also consider an \emph{extended model} which, in addition to the discrete probability part, also includes a continuum probability mass part. Then $\theta=(\theta_1,\theta_2,\dots)$, corresponding to the the discrete part of the distribution, only satisfies $\sum_\alpha \theta_\alpha \le 1$, where the remaining probability mass $\theta_0=1-\sum_\alpha \theta_\alpha$ belongs to the continuum part, the blob. We will derive the likelihood in this extended model and  define the PML $\hat{\theta}$ as the maximizer under the assumptions $ \theta_1\geq \theta_2\geq \ldots , \sum_{\alpha=1}^\infty \theta_{\alpha} \leq 1$. In Section \ref{Sec:Existence} we state the existence of the PML $\hat{\theta}$ in the extended model, and give the proof of this result in \ref{suppA}. Uniqueness is not known.

Both in the basic or extended model one can give a \emph{truncation level} $k=k_n$, and define $\widetilde \phi=(\theta_1,...,\theta_k)$ as well as $\phi_0=1-\sum_{\alpha=1}^k \theta_\alpha$. Such a truncated model we call a  \emph{sieved model}. Analogous to the standard PML one can write down a likelihood function for the sieved model and from this a PML, the so-called sieved PML. The introduction of the sieved PML (sPML) is novel and as discussed below is important for many applications.

The main theoretical results in the paper are almost sure consistency in an $L^{1}$-norm for the PML and sieved PML. In this connection the Hardy-Littlewood-Polya monotone rearrangement algorithm \cite{hardy:littlewood:polya:1952} is interesting for two reasons. The first reason is that the algorithm is prominent in our proof of the consistency result, since a naive estimator of the probability mass function can be seen as a monotone rearrangement of the empirical probability mass function. In the proof we need a certain contraction or non-expansivity property of the algorithm cf.\ \cite{anevski:fougeres:2007,lieb:loss:1996}. Another result is the almost sure rate of convergence which is almost of the order $n^{{-1/4}}$ for both the standard and sieved PML, which should be compared with the rate for the naive estimator, for which Jankowski and Wellner \cite{jankowski:wellner:2009} have obtained the rate $n^{-1/2}$, but then in distribution of norms, and furthermore for which we derive the almost sure rate in supnorm distance of order almost $n^{-1/2}$, cf. Section 3.   %\\

An important question is how to calculate the estimator. The main practical result is a stochastic approximation expectation maximisation (SA-EM) algorithm for the sieved estimator, where we use the EM algorithm to get a numeric approximation, treating the bijection $\chi$ as a latent variable; this is presented in   \ref{suppB}. In this algorithm, in the M step, assuming given $\chi$, we will use isotonic regression. We develop a modification of the standard PAVA algorithm for isotonic regression, cf.\ Robertson et al.\ \cite{robertson:wright:dykstra:1988}, to allow for lower bounds on the unknown frequencies, in \ref{suppC}.% the lower bounded isotonic regression estimator and the algorithm in detail in   C of \cite{anevski:gill:zohren:2013}.

The paper is organized as follows:  In Section \ref{Sec2} we introduce the model, the data that arise in this type of problem and the possible ways to estimate the probability mass function. In Section \ref{Sec:Existence} we state the existence result for the PML. In Section \ref{Sec3} we discuss consistency of the non-parametric maximal likelihood estimators: In Section \ref{Sec31} we will study an extended maximum likelihood estimator in the basic model, proving its consistency, and deriving rates  for the consistency result. In Section \ref{Sec32} we derive similar consistency results for the sieved estimator. In Section \ref{Sec33} we discuss the consistency results that we obtained in the previous two subsections and compare them with the results for the naive estimator obtained by Jankowski and Wellner \cite{jankowski:wellner:2009}. We conclude with a discussion in Section \ref{Sec5}. In   \ref{suppA} we prove existence of the PML. In   \ref{suppB} we present the SA-EM algorithm for computing the PML. In  \ref{suppC} we derive the MLE of a decreasing multinomial probability mass function bounded below by a known constant.

% ++++++++++++++++++++++++++++++++++++++++++++++++++++++++++++++++++++++++++
% ++++++++++++++++++++++++++++++++++++++++++++++++++++++++++++++++++++++++++
% ++++++++++++++++++++++++++++++++++++++++++++++++++++++++++++++++++++++++++
% ++++++++++++++++++++++++++++++++++++++++++++++++++++++++++++++++++++++++++
%\input section1

\section{The model, the data and the estimators} \label{Sec2}
\subsection{Introduction}

Imagine an area inhabited by a population of animals which can be classified by species. 
Which species live in the area (many of them previously unknown to science) is a priori unknown. 
Let  $\mathcal A$ denote the set of all possible species potentially living in the area.  
For instance, if animals are identified by their genetic code, then the species' names $\alpha$ are equivalence classes of DNA sequences. 
The set of all possible DNA sequences is effectively uncountably infinite, and for present purposes so is the set of equivalence classes, each equivalence class defining one potential species.  

Suppose that animals of species $\alpha\in\mathcal A$ form a fraction $\theta_\alpha \ge 0$ of the total population of animals. We assume that the probabilities $\theta_\alpha$ are unknown. The \emph{basic model} studied in this paper assumes that $\sum_{\alpha:\theta_\alpha>0} \theta_\alpha =1$ but we shall also study an \emph{extended model} in which it is allowed that (the discrete part of the distribution) $\sum_{\alpha:\theta_\alpha>0} \theta_\alpha < 1$. In either case, the set of species with positive probability is finite or at most countably infinite.

Imagine now an ecologist taking an i.i.d.\  random sample of $n$ animals, one at a time. The $j$th animal in the sample belongs to species $\alpha$ with probability $\theta_\alpha$. For each animal in turn, the ecologist can only determine whether it belongs to the same species as an earlier animal in the sample, or whether it is the first representative in his sample of a new species. Suppose he labels the different species observed in the sample by their number in order of discovery. His data can then be represented as a string of $n$ integers, where the $j$th integer equals $r$ if and only if it belongs to the $r$th different species observed in the sample in order of discovery. For instance, for $n=5$, the observed data could be the string $12231$ meaning that the first, second and fourth animals in the sample belonged to new species; the third and the fifth were each occurrences of a previously observed species, namely the same as that of the second and first animal in the sample respectively. 

%Since the sample is i.i.d., the data can be further reduced, by sufficiency, to the \emph{partition}, in the number-theoretic sense, of the integer $n$ which it induces. 
%This is the finite list of non-increasing integers $N=(N_1,N_2,\dots)$, with $\sum N_j=n$, where $N_j=m$ means that the $j$th most frequent species in the sample was observed exactly $m>0$ times (the definition does not depend on how ties are resolved). 
%For instance, the string $12231$ corresponds to the partition $N=(2,2,1)$ of the integer $5$, meaning that two species were each observed twice and one species was observed once; 2+2+1=5. It is convenient to append an infinite list of zero counts to $N$. In our example we then write $N=(2,2,1,0,0,\dots)$
\subsection{Estimation in the extended model}
Since we treat the $\alpha$ as unknown, the parameter $(\theta_\alpha:\alpha\in\mathcal A)$ is not identified. Since everything only depends on the ordered list of probabilities $\theta_\alpha$ it is convenient to change notation and from now on refer to species by their position in this ordering. If there is only a finite number of species of positive probability, then we will append to the list a countable number of possibly fictitious species each of probability zero. 
We redefine $\mathcal A =\mathbb N =\{1,2,\dots\}$ and redefine $\theta_\alpha$, where $\alpha$ is a positive integer, as the probability of the $\alpha$th most frequent species in the population. 
We define the \emph{deficit} $\theta_0=1-\sum_{\alpha\geq 1} \theta_\alpha$. 
In the basic model, $\theta_0=0$, in the extended model $\theta_0\ge0$. 

In the extended model, the deficit $\theta_0$ equals the probability, when we observe just one animal, that it belongs to one of those species which individually each have zero probability. 
Each such species can only be observed at most once in a sample of $n$ animals. 
The converse is not true: if an animal is observed only once in our sample, we do not know whether it belongs to a zero probability species or to a positive probability species.
 
We will discuss estimation in the extended model and in a truncated, or sieved, version of the extended model. 

Let $\aleph$ be the total number of species of positive probability. 
If $\aleph<\infty$, we take $\theta_\alpha=0$ for $\alpha >\aleph$. Thus from now on $\mathcal A=\mathbb N = \{1,2,...\}$, and $\theta=(\theta_1,\theta_2,...)$ where the $\theta_\alpha$,  the probability of occurrence of an animal belonging to the $\alpha$th most frequent species in the population, are nonnegative and nonincreasing and sum to 1. 

Since our random sample of $n$ animals is i.i.d., it can be further reduced, by sufficiency, to the \emph{partition}, in the number-theoretic sense, of the integer $n$ which it induces.  This is a list $N=(N_1,N_2,...)$ where $N_i\ge 0$ is the number of observed animals belonging to the $i$th most frequent species in the sample, $N_i\ge 0$, $N_1\ge N_2\ge...$, and $\sum_i N_i = n$.  The number $K$ of different species of animals observed in the sample, is finite: for some $K\ge 0$, $N_K>0$ and $N_i=0$ for $i > K$. In the number-theoretic sense of the word, $N$ (more precisely, the positive part of $N$, of length $K$) is a random \emph{partition} of the number $n$. For instance, the string $12231$ corresponds to the partition $N=(2,2,1)$ of the integer $5$, meaning that two species were each observed twice and one species was observed once; 2+2+1=5. It is convenient to append an infinite list of zero counts to $N$. In our example we then write $N=(2,2,1,0,0,\dots)$.

Both the data $N$ and unknown parameter $\theta$ are represented by infinite lists of nonincreasing nonnegative numbers, summing to $n$ and 1 respectively; the elements of $N$ are moreover integers. 
However there is no direct connection between the indices of the two lists. 
There exists a bijection $\chi$ from $\mathbb N$ (the species as ordered by the sample frequencies) to $\mathcal A$ (the species as ordered by population probabilities), defined by $\chi(i)=\alpha$ if and only if the $i$th most frequent species in the sample is the $\alpha$th most frequent species in the population. 
The bijection $\chi$ is random, and the essential feature of our model is that $\chi$  is not observed.

Let us use the same symbol $N$ to denote both the observed partition of sample size $n$ thought of as a random sequence, as well as the possible sample values thereof. 
After reduction by sufficiency, the sample space is the set of all possible partitions $N$ of the sample size $n$. 
Write $\textrm{P}^{(n,\theta)}$ for the corresponding (discrete) probability measure on the sample space when the underlying parameter is $\theta$. The basic model states that for any set $A$ of partitions of $n$
\begin{equation}
\textrm{P}^{(n,\theta)}(A)~=~\sum_{(N_1,N_2,...)\in A} {n \choose N_1~ N_2~ \dots} \sum_\chi \prod_i  \theta_{\chi(i)}^{N_i} .
\label{eq:Pdata}
\end{equation}
The likelihood function for $\theta$ based on the data $N$ is therefore
\begin{equation}
\textrm{lik}(\theta)~=~\sum_\chi \prod_i \theta_{\chi(i)}^{N_{i}}~=~\sum_\chi \prod_\alpha \theta_\alpha^{N_{\chi^{-1}(\alpha)}}.\label{eq:lik}
\end{equation}
The maximum likelihood estimator (MLE) of $\theta$ is defined as
\begin{equation}
\widehat\theta=\textrm{arg}\max_{\theta: \theta_1\geq \theta_2\geq \ldots, \sum_{\alpha=1}^\infty \theta_\alpha=1} \textrm{lik}(\theta). \label{eq:thetahat}
\end{equation}

It is interesting to note that the likelihood \eqref{eq:lik} can be interpreted as a matrix permanent of the nonnegative matrix $M_{ij} := \theta_i^{N_j}$. This relation enables one to use several techniques of approximate inference to evaluate the likelihood \cite{vontobel:2012,vontobel:2014}. We will not pursue this idea further here. This is mainly because we are interested in the extended model, where a relation to matrix permanents is more involved. 

Returning to the MLE, it is not clear that $\widehat \theta$ exists nor that it is unique.
In fact, it is easy to exhibit observed data $N$ for which it does not exist; for instance, with $n=2$, the partition $N=(1,1)$, see \ref{suppA} for the simple demonstration. For this reason we study instead the extended model MLE. Define the extended model MLE or the {\em Pattern Maximum Likelihood estimator} (PML) as 
\begin{equation}
\widehat{\theta}=\textrm{arg}\max_{\theta: \theta_1\geq \theta_2\geq \ldots , \sum_{\alpha=1}^\infty \theta_{\alpha} \leq 1} ~\sum_\chi \frac{n!}{N_0! \prod_{i\geq 1} N_i !} \; \theta_0^{N_0} \prod_{\alpha=1}^{\infty}   \theta_\alpha^{N_{\chi^{-1}(\alpha)}},\label{eq:extendedMLE}
\end{equation}
with $N_{0}=n-\sum_{\alpha=1}^{\infty} N_{\chi^{-1}(\alpha)}$ and $\theta_0=1-\sum_{\alpha\geq 1} \theta_\alpha$. The mappings $\chi:\mathbb N\to \{0,1,\dots,\infty\}$ satisfy that for every $\alpha\ge 1$ there exists exactly one $i$ such that $\chi(i)=\alpha$, and that $\chi(i)=0$ implies $N_i=0$ or $1$. Note that since the data ends in a block of 1's, with $N_0$ of them belonging to blob species, $\sum_{i\geq 0}N_i = n + N_0$. Furthermore $n! / N_0 ! \prod_{i: \chi(i) \ne 0} N_i ! = n! / N_0 ! \prod_{i > 0} N_i ! $, since $ 1! = 1$.  According to Theorem 1 in \cite{orlitsky:santhanam:viswanathan:zhang:2004}, it is true in this \emph{extended} model that a maximum likelihood estimator does exist; moreover they claim in Corollary 5 that the support of the PML (the number of indices for which $\hat{\theta}_\alpha$ is positive) is finite.  We prove that the PML $\hat{\theta}$ exists in Section \ref{Sec:Existence}, although the uniqueness is not known. The probability measure corresponding to a possibly defective probability $\phi$ is given by, for any set $A$ of partitions of $n$,
\begin{equation}
\textrm{P}^{(n,\phi)}(A)~=~\sum_{(N_1,N_2,...)\in A} \sum_\chi  \frac{n!}{N_0! \prod_{i\geq 1} N_i !} \;  \phi_0^{N_{0}} \prod_{\alpha=1}^{\infty}   \phi_\alpha^{N_{\chi^{-1}(\alpha)}}, \label{eq:extendedL}
\end{equation}
with $N_{0}=n-\sum_{\alpha=1}^{\infty} N_{\chi^{-1}(\alpha)}$ and $\phi_0=1-\sum_{\alpha\geq 1} \phi_\alpha$.

The underlying permutation of species generated by our finite sample of animals is not observed. 
Had it been observed, we would have access to full data counts $X=(X_\alpha:\alpha\in A)$. 
Here, $X_\alpha=N_{\chi^{-1}(\alpha)}$ is the number of occurrences of species $\alpha$ in the sample. 
This ``underlying data'' has the multinomial distribution with parameters $n$ and $\theta$. 

For any summable list of nonnegative numbers $a=(a_1,a_2,...)$, denote by $T(a)$ the \emph{monotone rearrangement map} which rewrites the components of $a$ in decreasing order. 
The relation between the actually observed $N$ and the underlying data $X$ is very simply $N=T(X)$.

To the underlying multinomial count vector $X$ we associate the empirical cumulative distribution function $F^{(n)}$ of the observed animals' true species label-numbers $\alpha$, defined by $F^{(n)}(x)=n^{-1}\sum_{\alpha\le x} X_{\alpha}$.
Alongside this we define the empirical probability mass function $f^{(n)}$, thought of as a vector or list rather than a function, $f^{(n)}_\alpha=X_\alpha /n=F^{(n)}(\alpha)-F^{(n)}(\alpha-1)$. Finally, we define 
$$
\widehat f^{(n)}~=~ N/n~=~T(f^{(n)})
$$ 
the \emph{naive estimator} of $\theta$. The two ways we have expressed it, show that it is simultaneously the \emph{ordered empirical} probability mass function of the \emph{underlying} data, as well as being a \emph{statistic} in the strict sense -- a function of the actually observed data $N$.

The naive estimator $\widehat f^{(n)}$ of $\theta$ is a random element on our sample space of random partitions. Our main tool in proving $L_1$ consistency of the PML $\widehat \theta$ will be finding an observable event $A$, i.e., a subspace of the set of all possible sample outcomes, which has large probability under $\textrm{P}^{n,\theta}$, where $\theta$ is the true value of the parameter, but small probability under $\textrm{P}^{n,\phi}$, for all $\phi$ outside of a small $L_1$ ball around $\theta$. This event $A$ will be defined in terms of $\widehat f^{(n)}$ and of the true parameter $\theta$; in fact, it will be the event that $\widehat f^{(n)}$ lies within a certain small $L_\infty$ ball around $\theta$. Since this true value of $\theta$ is fixed, even if unknown to the statistician, there is no problem in using its value in the definition of the event $A$.
\subsection{Sieved estimation in the extended model}

%We also study the maximum likelihood estimator in a \emph{truncated} (thus sieved) version of the \emph{extended} model, but we propose its use precisely when the \emph{basic} model is thought to be true.

In applications, maximization of the likelihood can be computationally very demanding. In the extended model, the parameter $\theta=(\theta_1,\theta_2,\dots)$ satisfies $\sum_\alpha \theta_\alpha \le 1$, and the total probability in the blob is $\theta_0=1-\sum_{\alpha\geq 1} \theta_\alpha$.  Whenever an animal is drawn from ``the blob'', it represents a new species in the sample, which is only observed exactly once. Thus when $\theta_0>0$ and $n$ is large, the observed partition $N$ tends to terminate in a long sequence of components $N_i$ all equal to $1$, many if not most of them -- in the long run, on average $\theta_0 n$ of them -- corresponding to species in the blob.

A possibly clever strategy for the \emph{basic} model would be to truncate the vector $\theta$ at some finite number of components. If however the true ordered probability mass function $\theta$ has a very slowly decreasing tail, truncation at too low a level might badly spoil the estimate. This possibility can be made less harmful by not truncating the original model, but truncating the extended model. Thus the parameter is taken to be $\widetilde\theta=(\theta_1,\dots,\theta_k)$ where $k<\infty$ and $\sum_{\alpha=1}^k \theta_\alpha\le 1$, and the probability deficit $\theta_0=1-\sum_{\alpha=1}^k \theta_\alpha$ is supposed to be spread  ``infinitely thinly'' over ``continuously many'' remaining species. 

These considerations lead to the idea of a sieved maximum likelihood estimator which we denote the {\em sieved PML estimator} (sPML), in which we maximize the probability of  the data over probability measures corresponding to a slightly different model from the true model, and indexed by a slightly different parameter: the model is both extended (to allow a blob) and truncated ($\theta$ has finite length). 

For given true parameter $\theta$ of basic or of extended model, and given truncation level $k=k_n$, define $\widetilde \theta=(\theta_1,...,\theta_k)$ and define $\theta_0=1-\sum_{\alpha=1}^k \theta_\alpha$. In general, $\widetilde \phi$ will denote a possibly defective probability mass function on $\{1,...,k\}$ where $\phi_1\ge\phi_2\ge ... \ge \phi_k$, and $\phi_0=1-\sum_{\alpha=1}^k \phi_\alpha$ will denote its deficit. Such parameters correspond to what we call the \emph{sieved model}.

Imagine the sieved model to be true. 
For any $i\in\mathbb N$, the species corresponding to the observed count $N_i\ge 0$ is either one of the species $\alpha=1,\dots,k$, or it is one of the species lumped together in the blob. The latter can only be the case if $N_i=1$ or $0$.
Different $i$ can both correspond to species in the blob, but cannot correspond to the same species in $1\le\alpha\le k$. 
We denote this mapping from $\mathbb N$ to $\{0,1,\dots,k\}$ by $\chi$. 
It cannot be a bijection, but every $1\le\alpha\le k$ does have a \emph{unique} inverse image. Moreover, $\chi(i)=0$ implies $N_i=1$ or $0$.
Apart from this it is arbitrary and not observed. 

Again we can imagine the full data which we would have had, if we had observed $\chi$. According to the sieved model there is an underlying $X=(X_0,X_1,...,X_k)$ which has the multinomial distribution with parameters $n$ and $(\phi_0,\widetilde\phi)$. To the ``proper part'' of $X$, that is to say, $(X_1,X_2,...,X_k)$, corresponds a partition of $X_+=\sum_{\alpha=1}^k X_\alpha$. Denote this partition by $N_+=(N_1,N_2, ... , N_J)$. Thus $J=\#\{1\le\alpha\le k: X_\alpha>0\}$ and $N_1\ge N_2\ge ... \ge N_J > 0$. Alongside these $X_+$ animals of $J\le k$ species from the set $\{1,\dots,k\}$, we also observed $X_0$ animals each of different species, where each of those species separately has probability 0, but all such species together have probability $\phi_0$. The observed data, finally, is the partition $N=(N_1,N_2, ... , N_J, 1,...,1)$ of $n$, in which we have appended exactly $X_0$ $1$'s to the partition $N_+$ of $X_+$.

Note that a number of the $N_i$ in the partition of $X_+$ can also equal $1$. In the observed data $N$ we cannot see how its block of $1$'s should be split between species inside and outside the blob. 

We can now write down the ``sieved likelihood'' and hence define the sPML estimator:
\begin{equation}
\textrm{lik}(\widetilde \phi)~=~\sum_\chi \frac{n!}{N_0! \prod_{i\geq 1} N_i !} \;\phi_0^{N_0} \prod_{\alpha=1}^k   \tilde{\phi}_\alpha^{N_{\chi^{-1}(\alpha)}},\label{eq:sievelik}
\end{equation}
\begin{equation}
\widehat\phi=\textrm{arg}\max_{\widetilde\phi: \widetilde\phi_1\geq \widetilde\phi_2\geq \ldots \geq\widetilde\phi_k, \sum_{\alpha=1}^k \widetilde\phi_\alpha\leq1} \textrm{lik}(\widetilde \phi).\label{eq:phihat}
\end{equation}
with $N_0=n-\sum_{\alpha=1}^k N_{\chi^{-1}(\alpha)}$ and $\phi_0=1-\sum_{\alpha=1}^k \tilde{\phi}_\alpha$. The mappings $\chi:\mathbb N\to \{0,1,\dots,k\}$ in the sum in (\ref{eq:sievelik}) have the properties that for every $1\le\alpha\le k$ there exists exactly one $i$ such that $\chi(i)=\alpha$, while $\chi(i)=0$ implies $N_i=0$ or $1$. It follows that the number of $i$ such that $N_i\ge 2$ cannot exceed $k$. 

Our strategy will again be to find an event $A$ such that $A$ has large probability under the true parameter but small probability under all parameters some distance from the truth. We do have to carefully distinguish between two different ``true'' probability measures: the law of the data within the sieved model, under the sieved parameter $\widetilde \theta$ corresponding to the truth, and the law of the data under the original, true model.

% ++++++++++++++++++++++++++++++++++++++++++++++++++++++++++++++++++++++++++
% ++++++++++++++++++++++++++++++++++++++++++++++++++++++++++++++++++++++++++
% ++++++++++++++++++++++++++++++++++++++++++++++++++++++++++++++++++++++++++
% ++++++++++++++++++++++++++++++++++++++++++++++++++++++++++++++++++++++++++
% ++++++++++++++++++++++++++++++++++++++++++++++++++++++++++++++++++++++++++

\section{Existence of the pattern maximum likelihood estimator}\label{Sec:Existence}

In this section we state an existence result for the PML estimator over an (extended) parameter space of ordered probability mass distributions in which we allow for a continuous part, the blob. We show existence by showing that this parameter space is compact, in an appropriate metric, and that the likelihood is a continuous functional with respect to this metric.

Recall that the extended parameter space $\Theta$ consists of sequences $\theta=(\theta_\alpha:\alpha\in \mathcal A)$ where $\mathcal A=\mathbb N=\{1,2,\dots\}$, and where $\theta_\alpha\ge 0$ for all $\alpha$, and moreover $\theta_1\ge\theta_2\ge \dots$  and $\sum_\alpha\theta_\alpha\le 1$.

We give  $\Theta$ the topology of pointwise convergence. Thus, for $\theta^{(m)}, \theta\in\Theta$, $\theta^{(m)}\to\theta$ as $m\to\infty$ if and only if $\theta^{(m)}_\alpha\to\theta_\alpha$ for all $\alpha$.

\bth $(i)$ Under the topology of pointwise convergence, the parameter space $\Theta$ is compact. $(ii)$ The functional $L: \Theta\mapsto {\mathbb R}_+$
defined by
\begin{eqnarray*}
   L(\theta)&=&~\sum_\chi  \frac{n!}{N_0! \prod_{i\geq 1} N_i !} \; \theta_0^{N_0} \prod_{\alpha=1}^{\infty}   \theta_\alpha^{N_{\chi^{-1}(\alpha)}}
\end{eqnarray*}
with $N_0=n-\sum_{\alpha=1}^{\infty} N_{\chi^{-1}(\alpha)}$, is continuous.

Thus the extended model pattern maximum likelihood estimator, defined in $(\ref{eq:extendedMLE})$, exists.
\eth
% ++++++++++++++++++++++++++++++++++++++++++++++++++++++++++++++++++++++++++
% ++++++++++++++++++++++++++++++++++++++++++++++++++++++++++++++++++++++++++
% ++++++++++++++++++++++++++++++++++++++++++++++++++++++++++++++++++++++++++
% ++++++++++++++++++++++++++++++++++++++++++++++++++++++++++++++++++++++++++
% ++++++++++++++++++++++++++++++++++++++++++++++++++++++++++++++++++++++++++
\section{Consistency results} \label{Sec3}
\subsection{Consistency for the PML estimator} \label{Sec31}

In this section we prove the consistency of the PML estimator in the extended model defined in $(\ref{eq:extendedMLE})$, based on a sample from the distribution $P$.  From our result of the previous section we know that there exists a PML. Uniqueness is not known; however our results below hold for any PML,  and in the sequel we let $\hat{\theta}$ denote {\em any} PML.

The idea of the proof is to first exhibit a sequence of events $A_n$ for which the $P^{n,\theta}$-probability is large (converges to 1 as $n\to \infty$), and such that for all probabilities $P^{n,\phi}$  such that $\phi$ is an $L_1$-distance $\delta$ away from $\theta$, the $P^{(n,\phi)}$-probability is small (goes to zero as $n\to\infty$). This is done in Lemma \ref{lem:1}.
 
As a consequence we show that the $P^{n,\theta}$-probability of $\{\frac{dP^{n,\phi}}{dP^{n,\theta}}>1$\}
is small (goes to zero as $n\to \infty$), by intersecting with $A_n$, for all $\phi$ that are $L_1$-distance more than $\delta$ away form $\theta$. On the other hand $\frac{dP^{n.\hat{\theta}}}{dP^{n,\theta}}>1$, if $\hat{\theta}$ is the ML estimator, for every ordered sample $(n_1,\ldots,n_k)$ with fixed $n=n_1+\ldots +n_k$. Finally we use an asymptotic formula for the number $p(n)$ of such $(n_1,\ldots,n_k)$, due to Ramanujan and Hardy, to make the argument uniform over every such sample, to show that $\hat{\theta}$ must be within $L_1$-distance of $\delta$ to $\theta$ with a large probability (that goes to one as $n\to\infty$), i.e. that $\hat{\theta}$ is weakly consistent. This is the content of Theorem \ref{thm:1}.

Using the bound established in Theorem \ref{thm:1}, we obtain almost sure consistency of $\hat{\theta}$, in Corollary \ref{cor:1}. Finally in Theorem 2 and Corollary 2, we derive rates of the almost sure convergence of the $L_1$ norm over classes of probability mass functions with tail conditions. 

Let $\theta$ be a fixed proper distribution. For $\delta>0$  arbitrary define the class of (possibly defective) probability mass functions ${\mathbb Q}_{\theta,\delta}=\{\phi: ||\phi-\theta||_1\geq \delta\}$, where $||\phi-\theta||_1=\sum_{i=1}^{\infty} |\phi_{i}-\theta_{i}|$. Note that $\phi$ is a possibly defective probability in the sense that $\sum_{i=1}^{\infty} \phi_{i}\leq1$, and note that in this case we use $(\ref{eq:extendedL})$ as the measure.

\blem\label{lem:1}
Let $f^{(n)}$ be the empirical probability mass function based on a sample $x_1,\ldots, x_n$ from some fixed decreasing probability mass function $\theta$, and $\hat{f}^{(n)}=T(f^{(n)})$. Then there is a finite $r=r(\delta,\theta)$ and $\epsilon=\delta/(8r)$ such that,
\begin{eqnarray*}
  { P}^{(n,\theta)}(\sup_{1\leq x\leq r} |\hat{f}^{(n)}_x-\theta_x|\leq {\epsilon})
&\geq &1- 2e^{-n\epsilon^2/2},\\
    \sup_{\phi\in {\mathbb Q}_{\theta ,\delta}}{P}^{(n,\phi)}(\sup_{1\leq x\leq r}|\hat{f}^{(n)}_x-\theta_x|\leq  {\epsilon}) &\leq & 2e^{-n\epsilon^2/2}.
\end{eqnarray*}
\elem
\prf
Let $\theta$ be fixed  and $\delta>0$ fixed but arbitrary, and choose an arbitrary $\phi \in {\mathbb Q}_{\theta,\delta}$. Since $\theta$ sums to one, there is an $r=r(\theta,\delta)$ such that $\sum_{i=r+1}^{\infty}\theta_i\leq \delta/4$. Then 
\begin{eqnarray}\label{eq:pq-diff}
     \sum_{i=1}^r|\theta_i-\phi_i|&\geq & \frac{\delta}{4}.
\end{eqnarray}
To show $(\ref{eq:pq-diff})$ note that either $\sum_{i=r+1}^\infty \phi_i$ is smaller or larger than $\delta/2$: $(i)$ Assume first that $\sum_{i=r+1}^\infty \phi_i\leq \delta/2$. Then 
\begin{eqnarray*}
     \delta&\leq& \sum_{i=1}^r|\theta_i-\phi_i|+\sum_{i=r+1}^\infty |\theta_i-\phi_i| \\
&\leq &\sum_{i=1}^r|\theta_i-\phi_i| +\sum_{i=r+1}^\infty \theta_i+\sum_{i=r+1}^\infty \phi_i\\
&\leq&\sum_{i=1}^r|\theta_i-\phi_i| +\frac{\delta}{4} +\frac{\delta}{2},
\end{eqnarray*}
which implies $(\ref{eq:pq-diff})$. $(ii)$ Assume instead that $\sum_{i=r+1}^\infty \phi_i> \delta/2$, and write the assumptions as $\sum_{i=1}^r \theta_i> 1-\delta/4$ and $\sum_{i=1}^r \phi_i= \sum_{i=1}^{\infty}\phi_{i}-  \sum_{i=r+1}^{\infty}\phi_{i} \leq 1-\delta/2.$  Then
\begin{eqnarray*}
     \sum_{i=1}^r |\theta_i-\phi_i|&\geq&  \sum_{i=1}^r (\theta_i-\phi_i) \\
&>& 1-\frac{\delta}{4}-1+\frac{\delta}{2}\\
&=&\frac{\delta}{4},
\end{eqnarray*}
which again implies $(\ref{eq:pq-diff})$. 

From $(\ref{eq:pq-diff})$ follows that for some $i\leq r$ we have 
\begin{eqnarray}\label{eq:pq-exist-i}
|\theta_i-\phi_i|&\geq& \frac{\delta}{4r}:=2\epsilon=2\epsilon(\delta,\theta).
\end{eqnarray}
Note that $r$, and thus also $\epsilon$ depends only on $\theta$, and not on $\phi$. 
%Note that the relation shows that given $\delta$, and $\epsilon>0$, by taking $r(P)$ large enough $(\ref{eq:pq-exist-i})$ holds (or if $r$ is fixed, by taking $\epsilon$ small enough).

Recall the Dvoretzky-Kiefer-Wolfowitz  (DKW) inequality \cite{dvoretzky:kiefer:wolfowitz:1956,massart:1990}; for every $\epsilon>0$
\begin{eqnarray}\label{eq:dkw:1}
     {\mathbb P}_{\theta}(\sup_{x\geq 0} |F^{(n)}(x)-F_{\theta}(x))|\geq {\epsilon}) &\leq &2 e^{-2n \epsilon^2},
\end{eqnarray}
where $F_{\theta}$ is the cumulative distribution function corresponding to $\theta$, and $F^{(n)}$ the empirical probability function based on i.i.d.\ data from $F_{\theta}$. Since $\{\sup_{x\geq 0}|F^{(n)}(x)-F_{\theta}(x)|\geq \epsilon\}\supset \{\sup_{x\geq 0}|f^{(n)}_x-\theta_x|\geq 2\epsilon\} \supset \{\sup_{x\geq 1}|f^{(n)}_x-\theta_x|\geq 2\epsilon\}$, with $f^{(n)}$ the empirical probability mass function corresponding to $F^{(n)}$, equation $(\ref{eq:dkw:1})$ implies
\begin{eqnarray}
  {P}^{n,\theta}(\sup_{x\geq 1}|f^{(n)}_x-\theta_x|\geq \epsilon)   &=&  {\mathbb P}_{\theta}(\sup_{x\geq 1}|f^{(n)}_x-\theta_x|\geq \epsilon)  \nonumber\\
     &\leq& 2e^{-n\epsilon^2/2 }. \label{eqbound}
\end{eqnarray}
Let $T$ be the monotone rearrangement map, cf.\ \cite{lieb:loss:1996}. Then the map $T$ is a contraction in the supnorm metric on ${\mathbb N}$ i.e.\, if $f,g$ are two functions ${\mathbb N}\to {\mathbb R}$ and $||f||_\infty=\sup_{k\geq 1}|f(k)|$ is the supnorm metric, then
$||T(f)-T(g)||_\infty\leq ||f-g||_\infty$, cf.\ \cite{anevski:fougeres:2007} (see also \cite{lieb:loss:1996}  for a proof
 of the contraction property for $L^p$-norms). Noting that $T(\theta)=\theta$ since $\theta$ is decreasing by assumption, and with $\hat{f}^{(n)}=T(f^{(n)})$, this implies that
\begin{eqnarray*}
    ||\hat{f}^{(n)}-\theta||_\infty&\leq & ||f^{(n)}-\theta||_\infty,
\end{eqnarray*}
so that $\{||\hat{f}^{(n)}-\theta||_\infty\geq \epsilon\}\subset\{|| {f}^{(n)}-\theta ||_\infty\geq \epsilon\}$, and thus by \eqref{eqbound}
\begin{eqnarray}\label{eq:PP-bound}
    {P}^{n,\theta}(\sup_{1\leq x\leq r} |\hat{f}^{(n)}_x-\theta_x|\geq {\epsilon}) &\leq &
 {P}^{n,\theta}(\sup_{x\geq 1} |\hat{f}^{(n)}_x-\theta_x|\geq {\epsilon}) \nonumber \\
&\leq & 2e^{-n\epsilon^2/2}.
\end{eqnarray}

For an analogue argument for a sample from the (possibly defective) distribution $\phi=(\phi_{1},\phi_{2},\ldots )$, we first append the mass point $\phi_{0}=1-\sum_{{x=1}}^{\infty}\phi_{x}$ to this vector to obtain a corresponding (proper) distribution function $F_{\phi}$. Using the corresponding cumulative empirical distribution $F^{(n)}$, and probability mass function $f^{(n)}$, and sorted such $\hat{f}^{(n)}=T(f^{(n)})$  we again have a contraction in the application of $T$, and going via the DKW inequality, we obtain (recall (\ref{eq:extendedL})),
\begin{eqnarray*}
 {P}^{n,\phi}(\sup_{1\leq x\leq r} |\hat{f}^{(n)}_x-\phi_x|\geq {\epsilon})
&\leq & 2e^{-n\epsilon^2/2},
\end{eqnarray*}
which is equivalent to
\begin{eqnarray}
        {P}^{n,\phi}(\sup_{1\leq x\leq r} |\hat{f}^{(n)}_x-\phi_x|< \epsilon)
&\geq & 1-2e^{-n\epsilon^2/2}. \label{eq:PQ-bound2}
\end{eqnarray}
Note that
\begin{eqnarray}
     && \{\sup_{1\leq x\leq r}|\hat{f}^{(n)}_x-\phi_x|< \epsilon\}\cap \{\exists i\leq r: |\theta_i-\phi_i|>2\epsilon\} \label{eq:6}\\
&&\subset  \{\exists i\leq r:|\hat{f}^{(n)}_i-\theta_i|> \epsilon\}=\{\sup_{1\leq x\leq r}| \hat{f}^{(n)}_x-\theta_x   |> \epsilon\}. \nonumber
\end{eqnarray}
Since the second event in $(\ref{eq:6})$ is deterministic, for any $\phi \in {\mathbb Q}_{\theta,\delta}$, and with an $\epsilon$ small enough (see $(\ref{eq:pq-exist-i})$), this together with equation $(\ref{eq:PQ-bound2})$ implies
\begin{eqnarray*}
    {P}^{n,\phi}(\sup_{1\leq x\leq r}|\hat{f}^{(n)}_x-\theta_x|>\epsilon)&\
\geq & {P}^{n,\phi}(\sup_{1\leq x\leq r}|\hat{f}^{(n)}_x-\phi_x|<\epsilon)\\
&\geq & 1-2e^{-n\epsilon^2/2}.
\end{eqnarray*}
Since $\phi\in{\mathbb Q}_{\theta,\delta}$ is arbitrary, the statement of the lemma follows.
\eop

We next derive the almost sure consistency of (any) extended maximum likelihood estimator $\hat{\theta}$. Recall the definitions of $P^{{n,\theta}}, P^{n,\phi}$ for proper and possibly defective distributions $\theta$ and $\phi$ in  $(\ref{eq:Pdata})$ and $(\ref{eq:extendedL})$, respectively.

\bth \label{thm:1}
Let $\hat{\theta}=\hat{\theta}^{(n)}$ be (any) extended maximum likelihood estimator. Then for any $\delta>0$ 
\[
{P}^{n,\theta}(||\hat{\theta}-\theta||_1>\delta) \leq \frac{1}{\sqrt{3}n} e^{\pi \sqrt{\frac{2n}{3}} -n\frac{\epsilon^2}{2}} (1+o(1))\quad \text{as} \quad n\to\infty
\]
where $\epsilon=\delta/(8r)$ and $r=r(\theta,\delta)$ such that $\sum_{i=r+1}^\infty \theta_i\leq \delta/4$.
\eth
%+++++++++++++++++++ Begin Proof Theorem 1
 \prf 
Now let ${\mathbb Q}_{\theta,\delta}$ be as in the statement of Lemma \ref{lem:1}. Then there is an $r$ such that the conclusion of the lemma holds, i.e.\, for each $n$ there is a set  
\begin{eqnarray*}
A&=&A_n=\{\sup_{1\leq x\leq r} |\hat{f}^{(n)}_x-\theta_x|\leq {\epsilon}\}
\end{eqnarray*}
such that
\begin{eqnarray*}
     { P}^{n,\theta}(A_n)&\geq & 1-2 e^{-n\epsilon^2/2},\\
     \sup_{\phi \in {\mathbb Q}_{\theta,\delta}}  { P}^{n,\phi}(A_n)&\leq& 2   e^{-n\epsilon^2/2}.
\end{eqnarray*}
For any $\phi\in {\mathbb Q}_{\phi,\delta}$, we can define the likelihood ratio $dP^{n,\phi}/dP^{n,\theta}$. Then for any  $\phi\in {\mathbb Q}_{\phi,\delta}$
\begin{eqnarray*}
     {P}^{n,\theta}\left(A_n\cap \left\{\frac{dP^{n,\phi}}{dP^{n,\theta}}\geq 1\right\}\right)&=& \int_{A_n\cap \left\{\frac{dP^{n,\phi}}{dP^{n,\theta}}\geq 1\right\}} dP^{n,\theta} \\
     &\leq & \int_{A_n} \frac{dP^{n,\phi}}{dP^{n,\theta}} \,dP^{n,\theta} \\
&=&{P}^{n,\phi}(A_n)\\
&\leq &2e^{-n\epsilon^2/2},
\end{eqnarray*}
which implies that
\begin{eqnarray*}
    {P}^{n,\theta}\left(\frac{dP^{n,\phi}}{dP^{n,\theta}}\geq 1\right)&=&{P}^{n,\theta}\left(A_n\cap \left\{\frac{dP^{n,\phi}}{dP^{n,\theta}}\geq 1\right\}\right)-{P}^{n,\theta}(A_n)\\
&&+{P}^{n,\theta}\left(A_n\cup \left\{\frac{dP^{n,\phi}}{dP^{n,\theta}}\geq 1\right\}\right)\\
&\leq & 2e^{-n\epsilon^2/2}-1+2e^{-n\epsilon^2/2}+1 \\
&=&4e^{-n\epsilon^2/2}.
\end{eqnarray*}
If $\hat{\theta}$ is a PML estimator then 
\begin{eqnarray*}
    \frac{dP^{n,\hat{\theta}}}{dP^{n,\theta}}&\geq &1.
\end{eqnarray*}
For a given $n=n_1+\ldots + n_k$ such that $n_1\geq \ldots \geq n_k>0$, (with $k$ varying), there is a finite number $p(n)$ of possibilities for the value of $(n_1,\ldots, n_{{k}})$. The number $p(n)$ is the partition function of $n$, for which we have the asymptotic formula
\begin{eqnarray*}
     p(n)&=&\frac{1}{4n \sqrt{3}}e^{\pi\sqrt{\frac{2n}{3}}}(1+o(1)),
\end{eqnarray*}
as $n\to \infty$, cf. \cite{ramanujan:hardy:1918}. For each possibility of $(n_1,\ldots, n_{{k}})$ there is a PML estimator (for each possibility we can choose one such) and we let ${\cal P}_n=\{\hat{\theta}^{(1)},\ldots, \hat{\theta}^{(p(n))}\}$ be the set of all such choices of PML estimators. Then
\begin{eqnarray*}
     {P}^{n,\theta}(\hat{\theta} \in {\mathbb Q}_{\theta,\delta})&=&\sum_{\phi\in {\cal P}_n \cap {\mathbb Q}_{\theta,\delta}} {P}^{n,\theta}(\hat{\theta} = \phi)\\
&\leq &\sum_{\phi\in {\cal P}_n \cap {\mathbb Q}_{\theta,\delta}} {P}^{n,\theta}\left( \frac{dP^{n,\phi}}{dP^{n,\theta}}\geq 1 \right) \\
&\leq &p(n) 4 e^{-n\epsilon^2/2},
\end{eqnarray*}
which ends the proof.
\eop
%+++++++++++++++++++ End Proof Theorem 1

That a $\hat{\theta}$ is consistent in probability is immediate from Theorem \ref{thm:1}, and in fact we have almost sure consistency:
 
\bcor\label{cor:1}
The sequence of maximum likelihood estimators $\hat{\theta}^{(n)}$ is strongly consistent in $L_1$-norm, i.e.\
\[
\lim_{n\to \infty}|| \hat{\theta}^{(n)}- \theta||_1 \stackrel{a.s.}{\to} 0
\]
as $n\to \infty$.
\ecor
\prf
This follows as a consequence  of the bound in Theorem \ref{thm:1}, by the characterization $X_n\stackrel{a.s.}{\to} 0 \Leftrightarrow \sum_{n=1}^{\infty}P(|X_n|>\delta)<\infty$ for all $\delta>0$, since
\begin{eqnarray*}
        \sum_{n=1}^{\infty}\frac{1}{\sqrt{3}n} e^{-\pi\sqrt{n}(\sqrt{n}\frac{\epsilon^2}{2}- \sqrt{\frac{2}{3}} })&<&\infty.
\end{eqnarray*}
\eop

The above results are for a fixed distribution $\theta$, and the rate depends, via $\epsilon$  on the distribution. The next Theorem and Corollary make the dependence explicit, and give a rate for the almost sure convergence as a function of the tail behaviour of the distribution.

\bth\label{thm:2}
Let $\epsilon_{0}>0$ be arbitrary and define 
\begin{eqnarray*}
{\Theta}_{\epsilon_{0}}=\{\theta: \forall \delta>0,\;  \exists r\leq \delta/\epsilon_{0}\mbox{ such that } \sum_{i=r+1}^{\infty}\theta_{i}<\delta/4\}.
\end{eqnarray*}
 Then, if $\theta\in {\Theta}_{\epsilon_{0}}$, 
\begin{eqnarray*}
         n^{\alpha}||\hat{\theta}^{(n)}-\theta|| &\stackrel{a.s.}{\to}& 0
\end{eqnarray*}
as $n\to \infty$, for any $\alpha<1/4$.
\eth
\prf
Let $\alpha>0$ be an arbitrary constant, to be determined below. From Theorem \ref{thm:1} we get
\begin{eqnarray}
{P}^{n,\theta}(n^{\alpha}||\hat{\theta}^{(n)}-\theta ||_1>\delta) &\leq &\frac{1}{\sqrt{3}n} e^{- n^{1/2}({n}^{1/2}\frac{\delta^2}{128\,  r^2 n^{2\alpha}}-\pi \sqrt{\frac{2}{3}} )}. \label{eq:r-bound}
\end{eqnarray} 
Since $\delta/r\geq \epsilon_{0}>0$ the right hand side of $(\ref{eq:r-bound})$ converges to zero, and is summable, if
\begin{eqnarray*}
{n^{-2\alpha+1/2}} &\to& \infty,
\end{eqnarray*}
as $n\to \infty$, which is true if  $\alpha< 1/4$.
\eop

\bcor\label{corr:tail}
Let $\Theta_{\kappa}=\{\theta:\theta_{x}=l(x)x^{{-\kappa}}\}$, for $\kappa>1$ fixed and with $l$ some function slowly varying at infinity. Then if $\theta\in \Theta_{{\kappa}}$ the conclusion of Theorem \ref{thm:2} holds.
\ecor
\prf
Assume that $\theta \in \Theta_{\kappa}$. Let $\epsilon_{0}>0$ be fixed, and let $\delta>0$ be fixed but arbitrary. Then for some $r$ we should have $\sum_{i=r+1}^{\infty}\theta_i<\delta/4$, which is equivalent to 
\begin{eqnarray*}
         r^{-\kappa+1} l_1(r) \leq \frac{\delta}{4} &\Leftrightarrow&  r\geq  (\frac{\delta}{4})^{1/(1-\kappa)} l_2(\delta), 
\end{eqnarray*}
when $\kappa >1$, where $l_1$ and $l_2$ are functions which vary slowly at infinity and zero respectively. It is possible to take $r$ such that $(\frac{\delta}{4})^{1/(1-\kappa)} l_2(\delta) \leq r <\delta/\epsilon_{0}$, thus $\theta\in \Theta_{\epsilon_{0}}$.
\eop

  % ++++++++++++++++++++++++++++++++++++++++++++++++++++++++++++++++++++++++++
% ++++++++++++++++++++++++++++++++++++++++++++++++++++++++++++++++++++++++++
% ++++++++++++++++++++++++++++++++++++++++++++++++++++++++++++++++++++++++++
% ++++++++++++++++++++++++++++++++++++++++++++++++++++++++++++++++++++++++++
% ++++++++++++++++++++++++++++++++++++++++++++++++++++++++++++++++++++++++++

%+++++++++++++++++++
\subsection{Consistency for the sPML estimator} \label{Sec32}
Let $k=k_n$ be a positive integer (truncation level) such that $k_n\to \infty$ when $n\to \infty$, and define the sieve 
\begin{eqnarray*}
    \widetilde{\Theta}_n&=&\{\widetilde \phi=(\phi_0,\phi_1,...,\phi_k) \mbox{ where }\phi_0=1-\sum_{\alpha=1}^k \phi_\alpha, \\
    && \mbox{and } \phi_{i}>\phi_{i+1},~ i=1,\ldots k-1\}.
\end{eqnarray*}  
Note that for each proper distribution $\phi\in \Theta_{\kappa}$ there is a corresponding sieved distribution $\tilde{\phi}\in  \widetilde{\Theta}_n$ with $\phi_0=\sum_{x=k_n+1} l(x)x^{-\kappa}\sim k_n^{-\kappa+1}$, if $\kappa>1$.

Assume the random vector  $X=(X_0,X_1,...,X_k)$, underlying our observations, has a multinomial distribution with parameters $n$ and $\widetilde\phi$. Define $J=\#\{\alpha\ge 1: X_\alpha>0\}$ and let $(N_1,N_2, ... , N_J)$ be a partition of $\sum_{\alpha=1}^k X_\alpha$, with $N_1\ge N_2\ge ... \ge N_J > 0$.  Then the observed data is the partition $(N_1,N_2, ... , N_J, 1,...,1,0,0,....)$ with $X_0\geq 0$ (unknown) number of $1$'s appended after the $J$'th position. Let $I=\sup\{i:N_i\geq 2\}$. We observe $I$, the number of species observed at least twice, and we observe $(J-I)+X_0$, the number of species which is only observed once. (We do not observe $J-I$ or $X_0$.) Note that the number of different species that we have observed frequency counts for is $J+X_0=\tilde{J}$, and that this number is known. We will let $k=k_n$ grow fast enough with $n$, so that always $\tilde{J}\leq k$.

Recall that $\chi:\{1,2,\ldots,\tilde{J} \}\to \{0,1,2,\ldots,k\}$ is a (random) map taking the $i$'th most frequently observed species to its position in the truncated list of species ordered by population frequency, such that all species above the $k$'th most common are grouped together in a ``zero category''. We assume that for every $\alpha$ such that $1\le\alpha\le k$ there is exactly one $1\leq i\leq \tilde{J}$ such that $\chi(i)=\alpha$. All other $i\in \{1,\ldots,\tilde{J}\}$ are mapped to the zero category.  This means that $\chi$ is injective on ${\cal I}=\chi^{-1}(\{1,\ldots,k\})$ and zero on its complement, so $\chi({\cal I}^c)=0$. Since $\tilde{J}\leq k$,  $\chi$ need not be surjective. The number $|{\cal I}|$ of observed species that are mapped to an $\alpha$ in $\{1,\ldots,k\}$ is random, although we do know that $|{\cal I}|\leq k$.

Define the sieved maximum likelihood estimator
\begin{eqnarray}
 \hat{\theta}_{(s)}^{(n)} =\argmax_{ \widetilde{\phi} \in \tilde{\Theta}_n} \sum_\chi \frac{n!}{N_0! \prod_{i\geq 1} N_i !} \phi_0^{N_0} \prod_{\alpha=1}^k   \phi_\alpha^{N_{\chi^{-1}(\alpha)}},\label{eq:sml-def}
\end{eqnarray}
with the sum running over all $\chi:\{1,2,\ldots, \tilde{J}\}\to \{0,1,\ldots,k\}$ such that $\chi$ is injective on a subset ${\cal I} \subset \{1,2,\ldots,\tilde{J}\}$, $\chi({\cal I})=\{1,\ldots,k\}$ and $\chi({\cal I}^c)=0$, and $N_0= n-\sum_{\alpha=1}^k N_{\chi^{-1}(\alpha)}$.

If $\chi$ and ${\cal I}$ are arbitrary but fixed  we define the ``estimator'' ${f}^{(n,\chi)}$ of a probability mass function on $\{0,1,\ldots,|{\cal I}|\}$ by
\begin{eqnarray}
    {f}^{(n,\chi)}(j)&=&\left\{ \begin{array}{ll}
\sum_{i\in {\cal I}^c} \frac{N_\chi(i)}{n},&\mbox{for }j=0\\
    T(\frac{N_{\chi(i)}}{n}:i\in {\cal I}),& \mbox{for }j\in \{1,\ldots,|{\cal I}|\}.
\end{array}
\right. \label{eq:Thateta}
\end{eqnarray}
This is not a proper estimator, since we can not calculate it only on the basis on our data $(N_1,N_2, ... , N_J, 1,...,1,0,0,....)$: The map $\chi$ and therefore the set ${\cal I}$ can not be determined from the sample.

For a given $\chi$, let $r_{\chi}$ be the restriction of a function $g$ on $\{1,2,\ldots\}$ to the set $\chi({\cal I})$. Define the map ${T}_{\chi}$ on the set of functions $g$ on $\{1,2,\ldots\}$ as the concatenation of the map $g \to \sum_{\alpha \in {\chi({\cal I})}^c} g_{\alpha}$, with the map composition of $T$ with $r_{\chi}$, so that 
\begin{eqnarray*}
T_{\chi}(g)&=&( \sum_{\alpha \in {\chi(\cal I})^c} g_{\alpha},T(r_{\chi}(g))).
\end{eqnarray*}
Then
\begin{eqnarray}
        T_{\chi}&:&\{\mbox{pmf  on }\{1,2,\ldots \}\}\mapsto \{\mbox{pmf on }\{0,1,\ldots,  |{\cal I}|\},\nonumber\\
        &&\mbox{ ordered on }\{1,\ldots,|{\cal I}|\}\},\label{eq:T-def}.
 \end{eqnarray}
If  $f^{(n)}$ is the empirical probability mass function, based on a sample $x_1,\ldots,x_n$ of $\phi$, cf.~Section 2, then       
\begin{eqnarray*} 
        {f}^{(n,\chi)}&=&{T}_{\chi}(f^{(n)}).
\end{eqnarray*}
Furthermore, for every $\chi$, the map ${T}_{\chi}$ in $(\ref{eq:T-def})$ is a contraction, with the two spaces of probability mass functions equipped with  the norms $|| \theta||=\sup_{ x\geq 1}|\theta_x|$ and $|| \theta||=\sup_{ 0\leq x\leq |{\cal I}|}|\theta_x|$, respectively. In particular
\begin{eqnarray}\label{eq:contraction}
    \sup_{0\leq x\leq |{\cal I}|} | {T}_{\chi}(f^{(n)})_x-{T}_{\chi}(\theta)_x| &\leq &\sup_{ x\geq 1 } |f^{(n)}_x-\theta_x|.
\end{eqnarray}

To show $(\ref{eq:contraction})$, note first that ${T}_{\chi}(\theta)=(\sum_{\alpha \in {\chi(\cal I})^c}\theta_{\alpha} ,\theta(\chi({\cal I})))$, since $\theta$ itself is sorted on $\chi({\cal I})$, and therefore $T_{\chi}(\theta)=\theta$ on ${\cal I}$. Furthermore ${f}^{(n)}$ is mapped to $(\sum_{\alpha\in \chi({\cal I})^c} f^{(n)}_{\alpha},T({f}^{(n)}(\chi({\cal I}))))$.

Therefore
\begin{eqnarray*}
    &&  \sup_{0\leq x \leq |{\cal I}|} | {T}_{\chi}({f}^{(n)})_x-{T}_{\chi}(\theta)_x| \\ 
    &=&\max(|\sum_{\alpha\in \chi({\cal I})^c} f^{(n)}_{\alpha}-\sum_{\alpha \in {\chi(\cal I})^c}\theta_{\alpha} |, \sup_{1\leq x\leq |{\cal I}|}| {T}(r_{\chi}(f^{(n)}))_x-{T}(r_{\chi}(\theta))_x|) \\
 &\leq &\max(|\sum_{\alpha\in \chi({\cal I})^c} f^{(n)}_{\alpha}-\sum_{\alpha \in {\chi(\cal I})^c}\theta_{\alpha}|, \sup_{x\in \chi({\cal I})}| {f}^{(n)}_x-{\theta}_x|) \\
 &\leq &\max(\sup_{x\in \chi({\cal I})^c}| {f}^{(n)}_x-{\theta}_x|), \sup_{x\in \chi({\cal I})}| {f}^{(n)}_x-{\theta}_x|) \\
&= &\sup_{ x\geq 1 } |{f}^{(n)}_x-{\theta}_x|,
\end{eqnarray*}
where the first inequality follows since the restriction of $T$ to any subset, and thus also to $\chi({\cal I})$, {\em is} a contraction, and the second inequality by the triangle inequality and since the $l^1$ norm on $\chi({\cal I})^c$ is bounded by the max-norm over $\chi({\cal I})^c$. This shows that $(\ref{eq:contraction})$ holds.

Define next the estimator ${\check{f}}^{(n)}$ of a probability mass function on the set $\{0,1,\ldots,I\}$, so on the blob together with the set of species observed at least twice, by 
\begin{eqnarray}
     \check{f}^{(n)}(j)&=&\left\{ \begin{array}{ll}
\sum_{i= I+1}^k \frac{N_i}{n},&\mbox{for }j=0\\
    \frac{N_j}{n},& \mbox{for }j\in\{1,\ldots,I\}.
\end{array}
\right. \label{eq:tildefn-def}
\end{eqnarray}
Note that this is a proper estimator. We extend this to an estimator on all of $\{0,\ldots, |{\cal I}|\}$ by defining $\check{f}^{(n)}(j)=0$ for $I<j\leq |{\cal I}|$.
 
We now have the following Lemma for the (extended) estimator $\check{f}^{(n)}$:
\blem\label{lem:extended}
Let $f_n$ be the empirical probability mass function based on a sample $x_1,\ldots, x_n$ from a fixed decreasing probability mass function $\theta$, and let $\check{f}^{(n)}$ be as defined in $(\ref{eq:tildefn-def})$. For $\delta>0$  arbitrary define the class of probability measures ${\mathbb Q}_{P,\delta}=\{Q: ||Q-P||_1\geq \delta\}$. Then there is a finite $r=r(\delta,P)$ and $\epsilon=\delta/(8r)$ such that,
\begin{eqnarray*}
  {P}^{n,\theta}(\sup_{1\leq x\leq r} |\check{f}_x^{(n)}-\theta_x|\leq {\epsilon})
&\geq &1- 2e^{-n(\epsilon-\frac{1}{n})^2/2},\\
    \sup_{\phi\in {\mathbb Q}_{\theta,\delta}}{P}^{n,\phi}(\sup_{1\leq x\leq r}|\check{f}_x^{(n)}-\theta_x|\leq  {\epsilon}) &\leq & 2e^{-n(\epsilon+\frac{1}{n})^2/2}.
\end{eqnarray*}
\elem
\prf
Let $\chi$ and $I$ be the fixed random elements that correspond to the given sample. Recall that $\chi$ is unknown and $I$ is known. From Lemma 1, there is an $r$ such that the conclusion of that Lemma holds. 

We first claim that
\begin{eqnarray*}
    \sup_{1\leq x\leq |{\cal I}|} |{f}^{(n,\chi)}_x-\check{f}^{(n)}_x|&\leq &\frac{1}{n}.
\end{eqnarray*}
To see this note first that ${f}^{(n,\chi)}$ and $\check{f}^{(n)}$ are identical on the set of species $\{1,\ldots,I\}$ that are observed at least twice. Since $ \check{f}^{(n)}$ is zero on $\{I+1,\ldots,|{\cal I}|\}$ it is enough to show that ${f}^{(n,\chi)}(j)\leq 1/n$ for $j\in \{I+1,\ldots,|{\cal I}|\}$. But this follows by the construction of ${f}^{(n,\chi)}$.

Therefore, with $||f||=\sup_{1\leq x\leq k}|f(x)|$ and recalling that $ |{\cal I}|\leq k$, we have $ ||\check{f}^{(n)}-\theta||\leq  \frac{1}{n}+||{f}^{(n,\chi)}-\theta||$
so that
\begin{eqnarray*}
     \{ ||{f}^{(n,\chi)}-\theta ||\leq \epsilon\}&\subset &\{||\check{f}^{(n)}-\theta||\leq \epsilon +\frac{1}{n}\},
\end{eqnarray*}
and from Lemma 1, with $n$ large enough that $1/n <\epsilon$,
\begin{eqnarray*}
  {P}^{(n,\theta)}(\sup_{1\leq x\leq r} |\check{f}^{(n)}_x-\theta_x|\leq {\epsilon})
&\geq &1- 2e^{-n(\epsilon-\frac{1}{n})^2/2}.
\end{eqnarray*}
Similarly 
\begin{eqnarray*}
     \{ ||\check{f}^{(n)}-\theta||\leq \epsilon\}&\subset &\{||{f}^{(n,\chi)}-\theta||\leq \epsilon +\frac{1}{n}\},
\end{eqnarray*}
so that from Lemma 1
\begin{eqnarray*}
    \sup_{\phi\in {\mathbb Q}_{\theta,\delta}}{P}^{(n,\phi)}(\sup_{1\leq x\leq r}|\check{f}^{(n)}_x-\theta_x|\leq  {\epsilon}) &\leq & 2e^{-n(\epsilon+\frac{1}{n})^2/2}.
\end{eqnarray*}
\eop

We need to get a bound on the total variation distance between the two measures $P^{n,\theta}$ and $P^{n,\tilde{\theta}}$ with $\theta$ a parameter and $\tilde{\theta}$ a sieved parameter. In order to get such a bound we need to make a coupling of the two measures. In particular the two random partitions $N,\tilde{N}$ of $n$ will be defined on the same probability space.

Therefore let $\theta=(\theta_1,\ldots,\theta_n)$ with $\theta_1\leq \theta_2 \leq\ldots \leq \theta_{k-1}\leq \theta_k \leq \theta_{k+1}\leq \ldots \leq \theta_n$ be the ordered set of probabilities. Note that the cut-off point defining the sieve is  $k=k_n$. The underlying full data is 
\begin{eqnarray*}
     (X_1,\ldots,X_n)&\sim& \mathrm{Multi}(n,\theta),
\end{eqnarray*}
where the $X_i$'s can be zeros and they need  not be ordered. Now let $X_0=\sum_{i=k+1}^n X_i$ and define the new underlying data $\tilde{X}=(X_0,X_1,\ldots,X_k)$. Then 
\begin{eqnarray*}
     \tilde{X}&\sim& \mathrm{Multi}(n,\check{\theta})
\end{eqnarray*}
where
\begin{eqnarray*}
       \check{\theta}&=&(\sum_{i=k+1}^n\theta_i,\tilde{\theta}),\\
       \tilde{\theta}&=&(\theta_1,\ldots,\theta_k).
\end{eqnarray*}
Now $N$ is the random partition of $n$, defined as the ordered $(X_1,\ldots,X_n)$, and $\tilde{N}$ is the random partition of $n$, defined by the ordered non-zero $X_1,\ldots,X_k$, to which we append a list of 1's of length $X_0$. Note that $N$ and $\tilde{N}$ are defined on the same probability space. Next  for any set $A$ of partitions on $n$ we define the two measures $P^{(n,\theta)},P^{(n,\tilde{\theta})}$ by
\begin{eqnarray*}
\textrm{P}^{(n,\theta)}(A)~=~\sum_{(N_1,N_2,...)\in A} {n \choose N_1~ N_2~ \dots} \sum_\chi \prod_{i=1}^n  \theta_{\chi(i)}^{N_i}, \\
\textrm{P}^{(n,\tilde{\theta})}(A)~=~\sum_{(\tilde{N}_1,\tilde{N}_2,...)\in A} {n \choose \tilde{N}_1~ \tilde{N}_2~ \dots} \sum_\chi \prod_i  \theta_{\chi(i)}^{\tilde{N}_i},
\end{eqnarray*}
in the case that $\theta$ is a proper distribution, and similarly if $\theta$ is a possibly defective distribution.  Note that $P^{(n,\theta)},P^{(n,\tilde{\theta})}$ have total mass one and thus are probability measures. There is another measure,  $\tilde{P}^{(n,\tilde{\theta})}$ say,  not necessarily a probability measure and connected to $\textrm{P}^{(n,\tilde{\theta})}$, that is defined by distributing the sorted nonzero values of $X_1,\ldots,X_k$ to different $\theta_i$'s and the value $X_0$ to the blob $\theta_0$. However, since we are only interested in when the measure $\textrm{P}^{(n,\theta)}$ differs from "the measure" generated by the partition $\tilde{N}$, it will not be of importance which of the two measures $P^{(n,\tilde{\theta})},\tilde{P}^{(n,\tilde{\theta})}$ we use, and as a matter of fact using a measure with total mass one simplifies the reasoning somewhat, therefore we will work with $P^{(n,\tilde{\theta})}$.

Now $P^{(n,\theta)}$ and $P^{(n,\tilde{\theta})}$ are the same if and only if all $X_{k+1},X_{k+2},\ldots,X_n$ are zero or one, and thus they differ on the set $\cup_{i=k+1}^n \{X_i\geq 2\}$. The probability, under $\theta$, of this is
\begin{eqnarray*}
      P_{\theta}(\cup_{i=k+1}^n \{X_i\geq 2\})&\leq& \sum_{i=k+1}^n P_{\theta}\{X_i\geq 2\} \\
      &\leq &  \sum_{i=k+1}^n \frac{E_{\theta}(X_i)}{2} %\\&=&
      =\frac{n}{2}  \sum_{i=k+1}^n \theta_i,
\end{eqnarray*}
by Markov's inequality.

\bth\label{thm:3} Let $\hat{\theta}_{(s)}^{(n)}$ be the sieved PML estimator defined in $(\ref{eq:sml-def})$. Assume  the sieve cut-off $k(n)$ satisfies  $\sum_{i=k(n)+1}^{n}\theta_i\leq C e^{-\beta n^{1/2+\nu}}(1+o(1))$, as $n\to \infty$, for some $\nu,\beta>0$. Then for any $\delta>0$ 
\begin{eqnarray*}
&&{ P}^{(n,\theta)}(||\hat{\theta}^{(n)}_{(s)}-\tilde{\theta} ||_1>\delta) \leq\\
&& \frac{1}{2\sqrt{3}n}e^{\pi \sqrt{\frac{2n}{3}}} (e^{-n{(\epsilon+\frac{1}{n})^2}/{2}}+e^{-n{(\epsilon-\frac{1}{n})^2}/{2}}+C e^{-\beta n^{1/2+\nu}}) (1+o(1))
\end{eqnarray*}
as $n\to \infty$, where $\epsilon=\delta/(8r)$ and $r=r(P,\delta)$ such that $\sum_{i=r+1}^\infty \theta_i\leq \delta/4$, and $||\tilde{\theta}-\tilde{\phi} ||_1=\sum_{i=1}^k|\tilde{\theta}_i-\tilde{\phi}_i|$.
\eth
\prf
Lemma 2 implies that there is a set 
\begin{eqnarray*}
 A_n&=&\{\sup_{1\leq x\leq k_n} |\check{f}^{(n)}_x-\theta_x|\leq {\epsilon}\}
\end{eqnarray*}
such that
\begin{eqnarray*}
     { P}^{n.\theta}(A_n)&\geq & 1-2 e^{-n(\epsilon-\frac{1}{n})^2/2},\\
     \sup_{\phi \in {\mathbb Q}_{\theta,\delta}}  { P}^{n,\phi}(A_n)&\leq& 2   e^{-n(\epsilon+\frac{1}{n})^2/2}.
\end{eqnarray*}

Furthermore, under the assumption of the cut-off level $k(n)$ we have that 
\begin{eqnarray*}
P^{n,\tilde{\theta}}(A)-P^{n,\theta}(A)&\leq &e^{-\beta n^{1/2+\nu}}(1+o(1))
\end{eqnarray*}
as $n\to \infty$, for any event $A$, and any sieved parameter $\tilde{\theta}$.

Let $\tilde{\theta}$ be a sieved parameter, derived from $\theta$. For any $\phi$, with corresponding sieved parameter $\tilde{\phi}$ we can define the likelihood ratio $dP^{n,\tilde{\phi}}/dP^{n,\tilde{\theta}}$. Let ${\mathbb Q}_{\tilde{\theta},\delta}=\{\tilde{\phi}:||\tilde{\phi}-\tilde{\theta}||_1 >\delta\}$.  Then since $\{||\theta-\phi ||_1>\delta \}\supset \{||\tilde{\theta}-\tilde{\phi}||_1>\delta\}$, we have that $\tilde{\phi}\in {\mathbb Q}_{\tilde{\theta},\delta} \Rightarrow \phi \in {\mathbb Q}_{{\theta},\delta}$. Therefore, for any  $\tilde{\phi}\in {\mathbb Q}_{\tilde\theta,\delta}$,  the corresponding ${\phi}\in {\mathbb Q}_{{\theta},\delta}$, and
\begin{eqnarray*}
    {P}^{n,{\theta}}\left(A_n\cap \left\{\frac{dP^{n,\tilde{\phi}}}{dP^{n,\tilde{\theta}}}\geq 1\right\}\right) -Ce^{-\beta n^{1/2+\nu}}&\leq &{P}^{n,\tilde{\theta}}\left(A_n\cap \left\{\frac{dP^{n,\tilde{\phi}}}{dP^{n,\tilde{\theta}}}\geq 1\right\}\right)\\
     &=& \int_{A_n\cap \left\{\frac{dP^{n,\tilde{\phi}}}{dP^{n,\tilde{\theta}}}\geq 1\right\}} dP^{n,\tilde{\theta}} \\
     &\leq & \int_{A_n} \frac{dP^{n,\tilde{\phi}}}{dP^{n,\tilde{\theta}}} \,dP^{n,\tilde{\theta}} \\
&=&{P}^{n,\tilde{\phi}}(A_n)\\
&=&{P}^{n,{\phi}}(A_n)+Ce^{-\beta n^{1/2+\nu}}\\
&\leq &2e^{-n(\epsilon+\frac{1}{n})^2/2}+Ce^{-\beta n^{1/2+\nu}},
\end{eqnarray*}
which implies that
\begin{eqnarray*}
    {P}^{n,{\theta}}\left(\frac{dP^{n,\tilde{\phi}}}{dP^{n,\tilde{\theta}}}\geq 1\right)&=&{P}^{n,{\theta}}\left(A_n\cap \left\{\frac{dP^{n,\tilde{\phi}}}{dP^{n,\tilde{\theta}}}\geq 1\right\}\right)-{P}^{n,{\theta}}(A_n)\\
&&+{P}^{n,{\theta}}\left(A_n\cup \left\{\frac{dP^{n,\tilde{\phi}}}{dP^{n,\tilde{\theta}}}\geq 1\right\}\right)\\
&\leq & 2e^{-n(\epsilon+\frac{1}{n})^2/2}+2Ce^{-\beta n^{1/2+\nu}}-1+2e^{-n(\epsilon-\frac{1}{n})^2/2}+1 \\
&=& 2e^{-n(\epsilon+\frac{1}{n})^2/2}+ 2e^{-n(\epsilon-\frac{1}{n})^2/2}+2Ce^{-\beta n^{1/2+\nu}}.
\end{eqnarray*}
If $\hat{{\theta}}_{(s)}^{(n)}$ is the sieved PML estimator then
\begin{eqnarray*}
    \frac{dP^{n,\hat{\theta}_{(s)}^{(n)}}}{dP^{n,\tilde{\theta}}}&\geq &1.
\end{eqnarray*}
For a given $n=n_1+\ldots + n_k$ such that $n_1\geq \ldots \geq n_k>0$, (with $k$ varying), there is a finite number $p(n)$ of possibilities for the value of $(n_1,\ldots, n_{{k}})$, for which the asymptotic formula
\begin{eqnarray*}
     p(n)&=&\frac{1}{4n \sqrt{3}}e^{\pi\sqrt{\frac{2n}{3}}}(1+o(1)),
\end{eqnarray*}
as $n\to \infty$, cf. \cite{ramanujan:hardy:1918}, holds. For each possibility of $(n_1,\ldots, n_{{k}})$ there is a sieved PML estimator and we let ${\cal P}_n=\{\hat{\theta}_{(s)}^{(n),(1)},\ldots, \hat{\theta}_{(s)}^{(n),(p(n))}\}$ be the set of all possible sieved PML estimators. Then 
\begin{eqnarray*}
     {P}^{n,{\theta}}(||\hat{\theta}_{(s)}^{(n)}-\tilde{\theta}||_1>\delta)&=&\sum_{\tilde{\phi}\in {\cal P}_n \cap {\mathbb Q}_{\tilde{\theta},\delta}} {P}^{n,{\theta}}(\hat{\theta}_{(s)}^{(n)} = \tilde{\phi})\\
&\leq &\sum_{\tilde{\phi} \in {\cal P}_n \cap {\mathbb Q}_{\tilde{\theta},\delta}} {P}^{n,{\theta}}\left( \frac{dP^{n,\tilde{\phi}}}{dP^{n,\tilde{\theta}}}\geq 1 \right) \\
&\leq &2p(n)(e^{-\frac{n}{2}(\epsilon-\frac{1}{n})^2}+e^{-\frac{n}{2}(\epsilon+\frac{1}{n})^2}+Ce^{-\beta n^{1/2+\nu}}).
\end{eqnarray*}
This ends the proof.
  \eop

The sieved PML estimator is strongly consistent:
\bcor\label{cor:2}
Under the assumption of Theorem  $\ref{thm:3}$, the sequence of sieved maximum likelihood estimators $\hat{\theta}_{(s)}^{(n)}$ is strongly consistent in $L_1$-norm, i.e.\
\[
|| \hat{\theta}_{(s)}^{(n)}- \tilde{\theta}||_1 \stackrel{a.s.}{\to} 0
\]
as $n\to \infty$.
\ecor
\prf
Follows from Theorem  $\ref{thm:3}$, analogously to Corollary $\ref{cor:2}$.
\eop
Note that if $\theta\in \Theta_{\kappa}$, so that $\theta_x=l(x)x^{-\kappa}$ with $l(x)$ a function slowly varying at infinity and $\kappa>1$, then the condition on the cut-off point is
  \begin{eqnarray*}
       Ce^{-\beta n^{1/2+\nu}}&\sim& \sum_{i=k(n)+1}^n \theta_i\sim \sum_{i=k(n)+1}^n i^{-\kappa}=k(n)^{-\kappa}\sum_{i=1}^{n-k(n)} i^{-\kappa} \\
       &\sim& k(n)^{-\kappa}(n-k(n))^{-\kappa+1}\\
      &\geq &k(n)^{-\kappa}n^{-\kappa+1},
  \end{eqnarray*}
  where the last inequality follows since $\kappa>1$ and $k(n)<n$. There is no way that we can have the condition of Theorem  $\ref{thm:3}$ satisfied if we only assume $\theta \in \Theta_{\kappa}$.

\bth \label{thm:4}
Let ${\Theta}_{\nu,\beta}=\{\theta: \theta_x=o(x^{\nu-1/2} e^{-\beta x^{\nu +1/2}}) \mbox{ as }x\to \infty\}$ for  $\nu>0,\beta>0$ fixed. Then, if $\theta\in {\Theta}_{\nu,\beta}$, 
\begin{eqnarray*}
         n^{\alpha}||\hat{\theta}^{(n)}_{(s)}-\tilde{\theta}|| &\stackrel{a.s.}{\to}& 0
\end{eqnarray*}
as $n\to \infty$, with $\alpha<1/4$.
\eth
\prf
Assume that $\theta \in \Theta_{\nu,\beta}$. Then the condition on exponentially decreasing tails in Theorem $\ref{thm:3}$ is satisfied. Furthermore, the condition $\forall \delta >0\; \exists r<\infty$ such that $\sum_{x=r}^{\infty} \theta_x <\delta/4$,   translates to 
\begin{eqnarray*}
\delta/4 \geq  e^{-\beta r^{1/2+\nu}} &\Leftrightarrow& r\geq \left( \frac{-\log \delta/4}{\beta}\right)^{2/(1+2\nu)}.
\end{eqnarray*}
The dominant part of the exponent in the right hand side of Theorem $\ref{thm:3}$ is then, replacing $\delta$ with $\delta/n^{\alpha}$ for an $\alpha$ to be chosen and with $\epsilon=\delta/8r$ and $r\sim (-\log \delta)^{2/(1+2\nu)}$,
\begin{eqnarray*}
      n^{1/2}-n \epsilon^2-2\epsilon-1/n&\sim & n^{1/2} -\frac{n^{1-2\alpha} \delta^2}{(-\log \delta)^{4/(1+2\nu)}}  - \frac{n^{-\alpha}\delta}{(-\log \delta)^{2/(1+2\nu)} }  \\
      &=&n^{1/2} -n^{1-2\alpha} c_1(\delta) -n^{-\alpha} c_2(\delta),
\end{eqnarray*}
which converges to $-\infty$ as $n\to \infty$ if $1-2\alpha >1/2$ and $\alpha> 0 $ i.e. if $0<\alpha< 1/4$.
Thus the rate is $n^{\alpha}$ for  any $\alpha< 1/4$.
\eop
%
  % ++++++++++++++++++++++++++++++++++++++++++++++++++++++++++++++++++++++++++
% ++++++++++++++++++++++++++++++++++++++++++++++++++++++++++++++++++++++++++
% ++++++++++++++++++++++++++++++++++++++++++++++++++++++++++++++++++++++++++
% ++++++++++++++++++++++++++++++++++++++++++++++++++++++++++++++++++++++++++
% ++++++++++++++++++++++++++++++++++++++++++++++++++++++++++++++++++++++++++
 \subsection{Comparison to the naive estimator} \label{Sec33}
 
An alternative to the non-para\-met\-ric maximum likelihood estimators, studied in the previous two subsections, is the naive estimator, consisting of estimating first the order relation from the data, and then given that estimate the population frequency by the observed population frequencies. 

We can obtain stronger results for the naive estimator than for the non-parametric maximum likelihood estimators.  In fact we can state almost sure supnorm convergence of the naive estimator with an almost parametric rate.

\blem \label{lem:naive}
Let $\hat{f}^{(n)}=T(f^{(n)})$ be the naive estimator. Then for any $\epsilon>0$ 
\[
{P}^{n,\theta}(||\hat{f}^{(n)}-\theta||_{\infty}>\epsilon) \leq 2e^{ -n{\epsilon^2}/{2}}\]
\elem
%+++++++++++++++++++ Begin Proof Theorem 1
 \prf 
 We argue similarly to the proof of Lemma \ref{lem:1}: Combining the Dvoretzky-Kiefer-Wolfowitz  inequality \begin{eqnarray*}\label{eq:dkw:2}
     {\mathbb P}_{\theta}(\sup_x |F^{(n)}(x)-F_{\theta}(x))|\geq {\epsilon}) &\leq &2 e^{-2n \epsilon^2},
\end{eqnarray*}
with $\{\sup_x|F^{(n)}(x)-F_{\theta}(x)|\geq \epsilon\}\supset \{\sup_x|f^{(n)}_x-\theta_x|\geq 2\epsilon\}$, we get
\begin{eqnarray*}
     {\mathbb P}_{\theta}(\sup_x|f^{(n)}_x-\theta_x|\geq \epsilon)&=&   {P}^{n,\theta}(\sup_x|f^{(n)}_x-\theta_x|\geq \epsilon) \\
     &\leq& 2e^{-n\epsilon^2/2 }.
\end{eqnarray*}
From the contraction property $||T(f)-T(g)||_{\infty}\leq ||f-g||_{\infty}$ of the monotone rearrangement map $T$ and since  $T(\theta)=\theta$, with $\hat{f}^{(n)}=T(f^{(n)})$, this implies that $\{||\hat{f}^{(n)}-\theta||_\infty\geq \epsilon\}\subset\{|| {f}^{(n)}-\theta ||_\infty\geq \epsilon\}$ and 
\begin{eqnarray*}
    {P}^{n,\theta}(\sup_{x} |\hat{f}^{(n)}_x-\theta_x|\geq {\epsilon}) &\leq &2e^{-n\epsilon^2/2}.
\end{eqnarray*}
\eop
Lemma \ref{lem:naive} implies consistency in probability, with rate $\alpha(n)=n^{{1/2}}(\log{n})^{-1/2}$, since  then $e^{-n\epsilon^{2}/2\alpha(n)^{2}}=e^{-\epsilon^{2}\log{n}/2}=n^{-\epsilon^{2}/2}$, which goes to zero, for every $\epsilon$. Almost sure consistency with rate $\alpha(n)=n^{1/2+\delta}$ holds, since $e^{-n\epsilon^{2}/2\alpha(n)^{2}}=e^{-n^{\delta}\epsilon^{2}/2}$ which is summable (in $n$). 

Thus we have the almost sure convergence and convergence in probability
\begin{eqnarray*}
     n^{1/2-\delta}  ||\hat{f}^{(n)}-\theta||_{\infty}  &\stackrel{a.s.}{\to}& 0,\\
       \frac{n^{1/2}}{\log{n}^{1/2}}  ||\hat{f}^{(n)}-\theta||_{\infty}  &\stackrel{P}{\to}& 0,
\end{eqnarray*}
for any $\delta>0$, as $n\to \infty$,

For the sieved model, recall the definition $(\ref{eq:tildefn-def})$ of the estimator ${\check{f}}^{(n)}$. Then similarly to the proof of Lemma \ref{lem:extended} we obtain the following result.
\blem
Let $f_n$ be the empirical probability mass function based on a sample $x_1,\ldots, x_n$ from a fixed decreasing probability mass function $\theta$, and let $\check{f}^{(n)}$ be as defined in $(\ref{eq:tildefn-def})$. Then, for any $\epsilon>0$,\begin{eqnarray*}
  {P}^{n,\theta}(||\check{f}^{(n)}-\theta||_{\infty}> {\epsilon})
&\leq &2e^{-n(\epsilon-\frac{1}{n})^2/2}. 
\end{eqnarray*} %\eop
\elem
As a consequence, this again give above rates in the two convergence modes.

 % ++++++++++++++++++++++++++++++++++++++++++++++++++++++++++++++++++++++++++
% ++++++++++++++++++++++++++++++++++++++++++++++++++++++++++++++++++++++++++
% ++++++++++++++++++++++++++++++++++++++++++++++++++++++++++++++++++++++++++
% ++++++++++++++++++++++++++++++++++++++++++++++++++++++++++++++++++++++++++
% ++++++++++++++++++++++++++++++++++++++++++++++++++++++++++++++++++++++++++
 
%\input section_computation
%
%

\section{Discussion}   \label{Sec5}

We discuss a non-parametric maximum likelihood estimator (PML) for a probability mass function with unknown labels, an estimator first introduced in the computer science literature by Orlitsky et al. \cite{orlitsky:sajama:santhanam:viswanathan:2004a} under the name of high profile estimator.  In Section \ref{Sec2}, we also introduced a sieved estimator which has a truncation level on the size of the probability vector. The existence of the PML estimator is proven in   \ref{suppA}.

The possibility of extending the model to include a continuous probability mass was already mentioned in \cite{orlitsky:sajama:santhanam:viswanathan:2004a}, however, it was not pursued further there. The introduction of a sieved estimator on the extended model is new and as we discuss below is important for many practical applications. 

In Section \ref{Sec3}, we proved strong consistency of ``the'' (actually any) PML (Theorem \ref{thm:1} and Corollary \ref{cor:1}) and sieved PML (Theorem \ref{thm:3} and Corollary \ref{cor:2}). The consistency of the PML was already claimed in \cite{orlitsky:sajama:santhanam:viswanathan:2005} without complete proof. The key ingredients to prove Theorem \ref{thm:1} and \ref{thm:3} are Lemma \ref{lem:1} and \ref{lem:extended} respectively. Both Lemmas use a novel strategy in proving  consistency of the MPL by finding an observable event $A$, which has large probability under $\textrm{P}^{n,\theta}$, where $\theta$ is the true value of the parameter, but small probability under $\textrm{P}^{n,\phi}$, for all $\phi$ outside of a small ball around $\theta$. Besides strong consistency we also determined the rate of convergence of the regular and sieved PML in Theorem \ref{thm:2} and \ref{thm:4} respectively, which in both cases is almost of the order $n^{-1/4}$. We conclude Section \ref{Sec3} by giving an comparison to the naive estimator by proving a result analogous to Lemma \ref{lem:1} and \ref{lem:extended} for the latter. 

\begin{remark}
The obtained almost sure rate of convergence for the PML is (almost) $n^{-1/4}$. It is not clear what the optimal almost sure rate is:  From the results of \cite{jankowski:wellner:2009} the rate of convergence for the naive estimator is $n^{-1/2}$; however this is the distributional rate of the $L_{p}$ norms. The best possible almost sure rate for this problem could be $n^{-1/2}$, and it could be slower. From our own results in Section \ref{Sec33} we get almost sure rates $n^{-1/2+\delta}$ for any $\delta>0$ for the naive estimator, which is faster than the rates for our estimator, it is however not clear if this is the optimal rate. Concerning our estimator, either the rate we obtain is the right rate for the PML which would mean that the PML is not optimal. Or else, the approach we use for deriving the rates is not the strongest possible, and in fact the rate for the PML is faster than $n^{-1/4}$ and (perhaps) equal to the optimal.

One should also note that the standard approach to deriving best rates for estimators is to use more sophisticated methods, for instance localization techniques. Our method consists of giving maximal inequalities for each PML and combining the derived bounds with a bound on the number of such PML's. This is a crude method and it is perhaps even surprising that we obtain consistency and rates at all.
\end{remark}

Another major result is the introduction of an algorithm to numerically compute the sieved PML. This is presented in   \ref{suppB} where the computation is based on the stochastic approximation of an expectation maximisation algorithm (SA-EM). In \cite{orlitsky:sajama:santhanam:viswanathan:2004b} a Monte Carlo Hastings expectation maximisation algorithm (MH-EM) of the standard PML was given. Our main advancement over this work is that we introduced the algorithm for the sieved estimator, and that we improved the statistical part of the EM algorithm by using the stochastic approximation.

Using the sieved estimator instead of the extended standard estimator can be an advantage when there are many unknown species with correspondingly small probabilities in the populations. Such situation appear for example in forensic DNA analysis. 

We illustrate this advantage on a small data example: Consider the partition 6=3+1+1+1, i.e.\ one species was observed three times and three species were observed once. The solution to the estimation problem is  intuitive and can be proven analytically \cite{pan2009}: One species, say 1, has probability $1/2$ and there is a continuous probability mass with a total probability $1/2$, i.e.\ based on the data, when sampling a new element, one expects to obtain 1 again in half of the cases or to observe a new species in the other half of the cases. To derive this estimator numerically one would have to use the extended model and the here presented algorithm. Using the algorithm for the standard model and a number of species of order of the sample size, a uniform distribution over all species apart from species 1, would give a too big probability to each element. Similar situations occur in real data problems, i.e.\ situations in which one would like to choose the species size of order of the sample size, but still account for a large number of rare species which have a very small probability which is comparable in size among the rare species.

%A possible complications with boundary solutions of the optimisation problem is:

\begin{remark}
For the SA-EM algorithm we note that, for a given finite value of $K$ we know that for a given data set a maximum likelihood estimate of $\boldsymbol\theta$ does exist. For each smaller value of $K$ there will typically correspond another, necessarily different, maximum likelihood estimate. All these estimates, one for each value of $K$ up to some maximum, correspond to fixed points of the EM algorithm when run with a larger still value of $K$. The SAEM algorithm therefore has many possible limits, corresponding to all values of $K$ not larger than the value corresponding to \emph{the} maximum likelihood estimate of $K$ for the given data-set and also not larger than the value of $K$ chosen in the implementation of the algorithm. These limits lie on the boundary of the parameter space. Once the procedure has got rather close to the boundary of the parameter-space, it is very difficult to move away again, since the size of potential steps is continuously being made smaller through the weights $\gamma$. Another troublesome part of the boundary of the parameter space corresponds to a sequence of probabilities $p_a$ which are all equal to one another. For large problems, once a long stretch of equal probabilities has arisen, this long segment is very resilient to change. Only very slowly can it get longer or shorter (at either end). 
\end{remark}
 
Therefore, in some cases unwanted results (i.e.\ local maxima of the optimisation problem) can be obtained when moving close to the boundary of the parameter space, i.e.\ when components of the probability vector become zero. In those cases, the numerical estimation can be improved by explicitly putting a lower bound on the allowed components of the probability vector. This means that in the M step of the EM algorithm one should change the isotonic regression to an isotonic regression of a probability mass function with a lower bound. It turns out that this problem has not been addressed in the literature, see however Balabdaoui et al.\ \cite{balabdaoui:rufibach:santambrogio:2009} for the related problem in isotonic regression of a regression function, see also van Eeden \cite{vaneeden:1957} and \cite[Theorem 2.1]{robertson:wright:dykstra:1988}. We have given a full solution to the lower bounded isotonic regression of a probability mass function in   \ref{suppC}.

\paragraph*{Acknowledgements} SZ is currently supported by Nokia Technologies, Lockheed Martin and the University of Oxford. Early states of this work were partially supported by FAPERJ, CNPq and PUC-Rio. SZ also thanks the Mathematical Institute at Leiden University for kind hospitality. DA's research has been partially supported by the Swedish Research Council, whose support is gratefully acknowledged.

%+++++++++++++++++++++++++++++++++++++++++++++++++++++++++++++++++++++++++
%   References
%+++++++++++++++++++++++++++++++++++++++++++++++++++++++++++++++++++++++++
%\section*{References}

\begin{supplement}
\slink[url]{http://arxiv.org/abs/1312.1200}
\sname{Appendix A}\label{suppA}
\stitle{Existence of the PML}
%\sdescription{In this Supplement we prove existence of the extended model ' maximum likelihood estimator}
\end{supplement}

\begin{supplement}
\slink[url]{http://arxiv.org/abs/1312.1200}
\sname{Appendix B}\label{suppB}
\stitle{Computation of the PML}
%\sdescription{In this Supplement we discuss an implementation of data, of the likelihood and the Stochastic Approximation EM algorithm (SAEM) used to calculate the NPMLE introduced in the previous section, in particular, the sieved model defined in \eqref{eq:phihat}.}
\end{supplement}

\begin{supplement}
\slink[url]{http://arxiv.org/abs/1312.1200}
\sname{Appendix C}\label{suppC}
\stitle{An algorithm for estimating a decreasing multinomial probability with lower bound}
%\sdescription{In this Supplement we present an algorithm for bounded isotonic regression and its prove convergence.}
\end{supplement}

\addcontentsline{toc}{section}{References}

%\bibliography{hiprofile}

 \appendix
 
\section{Existence of the extended model nonparametric maximum likelihood estimator}\label{App:1}

We first give a simple demonstration of non-existence of the MLE in the basic model. Thus define
\begin{eqnarray*}
   \widehat\theta=\textrm{arg}\max_{\theta: \theta_1\geq \theta_2\geq \ldots, \sum_{\alpha=1}^\infty \theta_\alpha=1}~\sum_\chi \prod_i \theta_{\chi(i)}^{N_{i}}.
 \end{eqnarray*}
Assume $n=2$ and the partition $N=(1,1)$ The data give a likelihood 
\begin{eqnarray*}
      \sum_{\chi} \prod_{i=1}^{2} \theta_{\chi(i)}^{N_i}=2(\theta_{1}\theta_{2}+\theta_{1}\theta_{3}+\ldots +\theta_{2}\theta_{3}+\theta_{2}\theta_{4}+\ldots).
\end{eqnarray*} 
We see first that there can only be a solution if all $\theta$'s are equal. In fact, writing the half likelihood as
\begin{eqnarray*}
      \theta_{1}\theta_{2}+\theta_{1}(1-(\theta_{1}+\theta_{2}))+\theta_{2}(1-(\theta_{1}+\theta_{2}))+R
\end{eqnarray*}
where $R$ contains all terms with only indices $3$ and higher, and differentiating w.r.t. $\theta_1$ we see that there is a maximum if and only if $\theta_{2} +1-2\theta_{1}-\theta_{2}-\theta_{2}=0$, i.e. if and only if $1-2\theta_{1}-\theta_{2}=0$. Since the likelihood is symmetric in the parameters, we get that there is a maximum if and only if $1-2\theta_{i}-\theta_{j}=0$ for all $i\neq j$, which is only possible if all $\theta_i=\theta_j$.  But if the cardinality $|{\cal A}|$ of the species names  is infinite, the restrictions $\theta_i=\theta_j, i\neq j, \sum_{\alpha \in {\cal A}}\theta_{\alpha}=1$, are not satisfied for any choice of parameters. Therefore there is no (ordered) $\theta$ that maximizes the likelihood in this case.  If $|{\cal A}|=:\aleph<\infty$ however then there is a solution which clearly is $\hat{\theta}=(1/\aleph,\ldots, 1/\aleph)\in [0,1]^{\aleph}$.

Next we prove that the PML in the extended model always exists, as stated in Theorem 1 in the main article.
Recall the definition of the (PML) as 
\begin{equation}
\widehat{\theta}=\textrm{arg}\max_{\theta: \theta_1\geq \theta_2\geq \ldots , \sum_{\alpha=1}^\infty \theta_{\alpha} \leq 1} ~\sum_\chi \frac{n!}{N_0! \prod_{i\geq 1} N_i !}  \theta_0^{N_{0}} \prod_{\alpha=1}^{\infty}   \theta_\alpha^{N_{\chi^{-1}(\alpha)}},\label{eq:extendedMLE}
\end{equation}
with $N_{0}=n-\sum_{\alpha=1}^{\infty} N_{\chi^{-1}(\alpha)}$ and with the mappings $\chi:\mathbb N\to \{0,1,\dots,\infty\}$ satisfying that for every $\alpha\ge 1$ there exists exactly one $i$ such that $\chi(i)=\alpha$, and that $\chi(i)=0$ implies $N_i=0$ or $1$. Recall also the definition of the measure $\textrm{P}^{(n,\phi)}$ for the possibly defective probability $\phi$: For any set $A$ of partitions of $n$,
\begin{equation}
\textrm{P}^{(n,\phi)}(A)~=~\sum_{(N_1,N_2,...)\in A} \sum_\chi  \frac{n!}{N_0! \prod_{i\geq 1} N_i !}  \theta_0^{N_{0}} \prod_{\alpha=1}^{\infty}   \theta_\alpha^{N_{\chi^{-1}(\alpha)}}, \label{eq:extendedL}
\end{equation}
with $N_{0}=n-\sum_{\alpha=1}^{\infty} N_{\chi^{-1}(\alpha)}$.

Recall that $\Theta$ is given the topology of pointwise convergence. We would like to note that Orlitsky et al.\ \cite{pan2009} suggested that the $\ell_2$-norm does the job in deriving existence, and they are (almost) right.

\prf (Theorem 1)

$(i)$ To see that $\Theta$ is compact, consider a sequence $\theta^{(m)}$. For given $\alpha$ the sequence of numbers $\theta_{\alpha}^{(m)}$ is bounded, hence contains a convergent subsequence. By a standard diagonalisation argument, we can extract from $\theta^{(m)}$ a subsequence for which each coordinate converges. 

$(ii)$ 
Suppose we take an iid sample of size $n$ of animals of different species labeled $1,2,3,\ldots$. The species have probabilities $\theta_1\ge \theta_2\ge \dots$ where $\sum_{k=1}^\infty \theta_k = 1 - \theta_0$. 
The index $k=1, 2, ...$ labels species in (decreasing) order of their probabilities; $k=0$ stands for a ``blob'' of very many different species each of very small probability. 
Two different animals each given the species label $k=0$ will always belong to different species.

Let the r.v. $S_i$ denote the species label of the $i$th animal in our sample, $i =1,\dots, n$; $S_i \in \{1, 2, ...\}\cup\{0\}$, and note that $S_1,\ldots,S_n$ are i.i.d. r.v.'s. Note also that since the (theoretical) species labels are not observed, $S_1,\ldots,S_n$ are not statistics, they are however  random variables. When we have obtained our sample we can determine for any two elements of the sample whether they belong to the same species or not. This determines a random equivalence relation on the numbers $\{1, 2, \dots, n\}$, which we shall denote by $\sim$: $i\sim j$ if and only if $i=j$ or $i\ne j$ and $S_i = S_j \ne 0$. 

We will introduce a second random equivalence relation denoted by $\sim_K$: $i\sim_K j$ if and only if $i=j$ or $i\ne j$ and $S_i = S_j \in \{1,\dots,K\}$. These equivalence relations determine partitions ${\mathcal P}_n$ and ${\mathcal P}_n^K$ of the set $\{1, 2, \dots, n\}$ into equivalence classes, so e.g. ${\mathcal P}=\{H_1,\ldots, H_{\tilde{n}}\}$, with 
\begin{eqnarray*}
      \{1,\dots,n\}&=&\cup_{j=1}^{\tilde{n}} H_j,
\end{eqnarray*}
with $H_i\cap H_j=\emptyset$ if $i\neq j$. Note that the equivalence relation $\sim_K$ is \emph{stricter} than the equivalence relation $\sim $, in the sense that $i\sim_K j \Rightarrow i\sim j$, which implies that the partition ${\mathcal P_n}^K$ generated by $\sim_K$ is \emph{finer} than the partition ${\mathcal P}_n$ generated by  $\sim $, i.e. an equivalence set in $ {\mathcal P}_n$ is a union of equivalence sets of ${\mathcal P}_n^{K}$. 

The sizes of the equivalence classes determine partitions $n=|H_1|+\ldots +|H_{\tilde{n}}|$, in the number theoretic sense, of the number $n$.  Let $\Pi_n$ denote the random partition of the number $n$ generated by $\sim$ and $\Pi_n^K$ that generated by $\sim_K$. Denote by $\pi_n$ and $\pi_n^K$ possible realisations of both. Denote by $P_\theta$ the probability measure induced by $\theta = (\theta_1,\theta_2, \dots)$. Now given a parameter vector $\theta$ define $\theta^K = (\theta_1, \dots, \theta_K, 0, 0, \dots)$.  All species with label larger than $K$ have been merged with the blob. 

Define the event
\begin{eqnarray*}
A_{n,K} &=& \cup_{1\le i < j \le n}(\{S_i = S_j \}\cap\{S_i> K\}\cap\{S_j >K\}),
\end{eqnarray*} 
of at least two animals in the sample belong to the same species and have a species label larger than $K$. The complement is
\begin{eqnarray*}
   A_{n,K}^c&=&\cap_{1\le i < j \le n}(\{S_i \neq S_j \}\cup\{S_i\leq K\}\cup \{S_j\leq K\})
\end{eqnarray*}
i.e. the event that for every pair of animals no two are from the same species or at least one of the pair of animals has a label smaller than or equal to $K$. Note that 
\begin{eqnarray}
P_\theta(A_{n,K})     &\le &\frac 12 n(n-1) \sum_{i=K+1}^{\infty}\theta_{j}^{{2}}\nonumber\\
   &\leq &\frac 12 n(n-1) \theta_{K+1}\sum_{i=K+1}^{\infty}\theta_{j}\nonumber \\
  &\leq &\frac 12 n(n-1) \theta_{K+1},
 \label{eq:ProbAnk-bound}
\end{eqnarray}
where the first inequality follows by Boole's inequality and since the $S_i$ are i.i.d., and the second since $\theta_k \le \theta_{K+1}$ for $k\ge K+1$.
%\footnote{It does NOT hold that on the event $A_{n,K}^{\textrm c}$, also $\Pi_n = \Pi_n^K$. To see this: On the event $A_{n,K}^c$, for every pair $(i,j)$ either 1) $S_i\neq S_j$ in which case both of $i\sim j$ and $i\sim_K j$ are NOT satisfied, and $i,j$ are in different partition sets in both ${\mathcal P}_n$ and ${\mathcal P}_n^K$ or 2) $S_i=S_j$ and for instance $S_i\leq K, S_j>K$. In that second 2) case $i\sim j$ so $i,j$ are both in the same partition set in ${\mathcal P}_n$, but $i\sim_K j$ does NOT hold, so $i,j$ are not in the same partition set in ${\mathcal P}_n^K$. Thus $\Pi_n = \Pi_n^K$}

We have that $\Pi_n=\Pi_n^K$ on $A_{n,K}^c$. In fact, on $A_{n,K}^c$, let $1\leq i<j\leq n$ be fixed but arbitrary. Then, if $S_i\neq S_j$ both of $i\sim j$ and $i\sim_K j$ are violated so then $i,j$ are not in the same partition in ${\mathcal P}_n$ nor in ${\mathcal P}_n^K$. If instead $S_i=S_j$ then we must have that $S_j\leq K$ and $S_j\leq K$, and then if $S_i=S_j>0$ both $i\sim j$ and $i\sim_K j$ are satisfied so then $i,j$ are in the same partition in both ${\mathcal P}_n$ and ${\mathcal P}_n^K$, and if $S_i=S_j=0$ neither of $i\sim j, i\sim_K j$ are satisfied and then $i,j$ are not in the same partition in ${\mathcal P}_n$ nor in ${\mathcal P}_n^K$. Since this holds for every $i<j$, and since $\Pi_n$ and $\Pi_n^K$ are counting the sizes of the partitions in ${\mathcal P}_n$ and ${\mathcal P}_n^K$, we have shown that on $A_{n,K}^c$ the two partitions $\Pi_n$ and $\Pi_n^K$ of $n$, coincide. It therefore follows that for any given partition $\pi_n$ of the number $n$ 
\begin{eqnarray}
 P_\theta(\Pi_n = \pi_n) &=& P_\theta(\{\Pi_n = \pi_n\}\cap A_{n,K}^{\textrm c}) + P_\theta(\{\Pi_n = \pi_n\}\cap A_{n,K}) \nonumber \\
&\le&  P_\theta(\{\Pi_n^K = \pi_n\}\cap A_{n,K}^{\textrm c}) +P_\theta(A_{n,K}) \nonumber \\
& \le& P_\theta(\{\Pi_n^K = \pi_n\}) +P_\theta(A_{n,K}) \nonumber \\
&\le &P_{\theta^K}(\{\Pi_n = \pi_n\}) +\frac 12 n(n-1) \theta_{K+1}. \label{eq:pthetamanip}
\end{eqnarray}

%Now, since 
%\begin{eqnarray*}
%     A_{n,K}^c&=&\cup_{1\leq i<j\leq n} A_{n,K}^c\cap (\{S_i>K\}\cup \{S_j\leq K\})
%\end{eqnarray*}
%is a partition of $A_{n,K}^c$, the first term in $(\ref{eq:pthetamanip})$ can be bounded as
%\begin{eqnarray*}
%   &&P_\theta(\{\Pi^K_n = \pi_n\}\cap A_{n,K}^{\textrm c})\\
%&=&\sum_{1\leq i<j\leq n}\{ P_\theta(\{\Pi^K_n = \pi_n\}\cap A_{n,K}^{\textrm c}|S_i>K,S_j\leq K) P_{\theta}(S_i>K,S_j\leq K) \\
%&&+P_\theta(\{\Pi^K_n = \pi_n\}\cap A_{n,K}^{\textrm c},S_i\leq K,S_j\leq K)\}\\
%&=&\sum_{1\leq i<j\leq n}\{ P_\theta(\{\Pi^K_n = \pi_n\}\cap A_{n,K}^{\textrm c}|S_i>K,S_j\leq K) P_{\theta}(S_i>K)P_{\theta}(S_j\leq K) \\
%&&+P_\theta(\{\Pi^K_n = \pi_n\}\cap A_{n,K}^{\textrm c},S_i\leq K,S_j\leq K)\}
%\end{eqnarray*}
%where the first inequality holds since, for the first term, the probability distribution of $\Pi_n$ under $\theta^K$ is the same as the probability distribution of $\Pi_n^K$ under $\theta$ and for the second term the probability of an intersection is smaller than for any of the factors, and the second inequality follows by $(\ref{eq:ProbAnk-bound})$. }

Let $\theta^{(m)}$ be a sequence of parameter vectors converging coordinatewise to $\theta$ as $m \to \infty$, and let $\pi_n$ be a fixed partition of the number $n$. We want to prove that  $P_{\theta^{(m)}}(\Pi_n = \pi_n)\to P_\theta(\Pi_n = \pi_n)$ as $m\to\infty$.  

Let $\delta > 0$. Then there is finite $K=K(\delta)$ such that
\begin{eqnarray}
 \frac 12 n (n-1) \theta_{K+1} \le \frac 12 \delta. \label{eq:theta-bound}
\end{eqnarray}
Furthermore there is a finite $M=M(\delta,\theta_{K+1})$ such that if $m>M$, we have that  
\begin{eqnarray*}
 \frac 12 n (n-1) \theta_{K+1}^{(m)} &\le&  \frac{\delta}{2}. 
\end{eqnarray*}
Therefore, using the inequality$(\ref{eq:pthetamanip})$ with $\theta$ replaced by $\theta^{(m)}$,
\begin{eqnarray*}
    P_{\theta^{(m)}}(\Pi_n = \pi_n) &\leq & P_{\theta^{(m)K}}(\{\Pi_n = \pi_n\}) +\frac{\delta}{2}.
\end{eqnarray*}
Now $P_{\theta^{(m)K}}(\{\Pi_n = \pi_n\}) \to P_{\theta^{K}}(\{\Pi_n = \pi_n\})$ as $m\to \infty$, since $\theta^{(m)K}$ only contains finitely many non-zero coordinates. This implies that 
\begin{eqnarray}
  \limsup_{m\to \infty} P_{\theta^{(m)}}(\Pi_n = \pi_n) &\leq& P_{\theta^{K}}(\{\Pi_n = \pi_n\}) + \frac{\delta}{2}. \label{eq:limsupbound}
\end{eqnarray}

Next, we get
\begin{eqnarray}
   P_{\theta}(\{\Pi_n=\pi_n\})&\ge& P_{\theta}(\{\Pi_n=\pi_n\}\cap A_{n,K}^c) \nonumber\\
&= &P_{\theta}(\{\Pi_n^K=\pi_n\}\cap A_{n,K}^c) \nonumber \\
&= &P_{\theta}(\{\Pi_n^K=\pi_n\})-P_{\theta}(\{\Pi_n^K=\pi_n\}\cap A_{n,K})\nonumber \\
&\ge &P_{\theta^K}(\{\Pi_n=\pi_n\})-\frac{1}{2}n (n-1)\theta_K , \label{eq:pthetaeq2}
\end{eqnarray}
where the first equality follows since $\Pi_n^K=\Pi_n$ on $A_{n,K}^c$, and last inequality follows from $(\ref{eq:ProbAnk-bound})$.  Finally  $(\ref{eq:limsupbound})$, $(\ref{eq:pthetaeq2})$ and $(\ref{eq:theta-bound})$ imply that
\begin{eqnarray*}
     \limsup_{m\to \infty} P_{\theta^{(m)}}(\Pi_n = \pi_n) &\leq& P_{\theta}(\{\Pi_n = \pi_n\}) +\delta.
\end{eqnarray*}

To show a lower bound for the liminf, use of $(\ref{eq:pthetaeq2})$ with $\theta^{K}$ replaced by $\theta^{(m)K}$, and noting that $\theta_{K+1}^K=0$, gives
\begin{eqnarray}
    P_{\theta^{(m)K}}(\Pi_n=\pi_n)&\le& P_{\theta^{(m)}}(\Pi_n=\pi_n).\label{eq:thetaK-bound}
\end{eqnarray}
Thus
\begin{eqnarray*}
   P_{\theta^{(m)}}(\Pi_n=\pi_n)&\ge& P_{\theta^{(m)K}}(\Pi_n=\pi_n)\\
&\to& P_{\theta^K}(\Pi_n=\pi_n)\\
&\geq & P_{\theta}(\Pi_n=\pi_n)-\frac 12 \delta,
\end{eqnarray*}
where the first inequality holds by $(\ref{eq:thetaK-bound})$, then the limit (which is a liminf) is taken as $m\to \infty$ and the last inequality follows by $(\ref{eq:pthetamanip})$ and $(\ref{eq:theta-bound})$.

Thus lim sup and lim inf of $P_{\theta^{(m)}}(\Pi_n = \pi_n)$ are within $\delta$ of $P_\theta(\Pi_n = \pi_n)$.  
Since $\delta > 0$ was arbitrary it follows that $P_\theta(\Pi_n = \pi_n)$ is the limit as $m \to \infty$ of $P_{\theta^{(m)}}(\Pi_n = \pi_n)$.

\eop

\section{Computation of the nonparametric maximum likelihood estimator} \label{App:2}

In this appendix we discuss an implementation of data, of the likelihood and the Stochastic Approximation EM algorithm (SAEM) used to calculate the NPMLE introduced in the previous section, in particular, the sieved model defined in Equation (7) of the main article. %\eqref{eq:phihat}. 

\subsection{The sample}

After reduction by sufficiency, the data can be represented by the \emph{partition} of the sample-size $T$, in the number theoretic sense: A partition of $T$ is a non-increasing sequence of positive integers adding to up $T$, e.g., $T=7=3+2+1+1$. The number of different integers appearing in the partition can be much smaller than the length of the partition itself, and often a more compact representation of $T$ consists of two equal length sequences of positive integers $n_1<  \dots < n_J$ and $r_1,\dots,r_J$ where $n_j$ are the distinct numbers occuring in the partition, ordered, $r_j$ are the number of repetitions of $n_j$ and $J$ is the number of distinct numbers occurring in the partition. Write $\mathbf r = (r_j)_{1\le j\le J}$ and $\mathbf n=(n_j)_{1\le j\le J}$. In the above example $J=3$, $\mathbf n = (1,2,3)$, and $\mathbf r = (2,1,1)$.

\bass \label{ass:1}
Assume that $n_1=1$, i.e. there exist singletons in the sample, and  $J\ge 2$, i.e. the sample contains non-singletons.
\eass

Assumption \ref{ass:1} is typically satisfied in practice; in the sequel we assume this to hold.

\subsection{The population}

We will use indices $a$, $b$, etc.\ to denote (non-blob) \emph{population species}, identified by position when ordered by decreasing probability. Different blob species are merged into one group and assigned the index $0$. We suppose the population consists of a finite number $K$ of species of positive probability $p_1\ge p_2 \ge\dots\ge p_K > 0$ and a blob of uncountably many species each of zero probability, but together of positive probability $p_0=1-\sum_{a=1}^K p_a >0$. The population species $a$, $b$ etc.\  are therefore integers between $0$ and $K$ where $0$ indicates a blob species and $1$ to $K$ a non-blob species. 

In some situations one can be interested in the case $K=0$ but this special case is easy to study separately, so we will assume in the sequel $K\ge 1$.

\subsection{The likelihoods}

The ``missing data'' consists of the identification of each non-blob population species \emph{either} with an index $1\le j \le J$ to indicate that this species was indeed observed in the sample, and was one of the $r_j$ species observed exactly $n_j$ times, \emph{or} with some kind of marker, we will use the index $0$ for this purpose, to indicate that this species was not observed at all. 

Under Assumption \ref{ass:1} the number of singletons  $r_1$ in the sample is positive, and $J\ge 2$, so the sample contains both singletons and non-singletons.
Then, the missing data can be represented by a function $\boldsymbol\psi:\{1,\dots,K\}\to \{0,1,\dots,J\}$, which satisfies the two constraints
\begin{itemize}
\item[\bf C1:] $\sum_{a=1}^K 1\{\boldsymbol\psi(a)=j\}=r_j$, for each $j>1$,
\item[\bf C2:] $\sum_{a=1}^K1\{\boldsymbol\psi(a)=1\}\le r_1$.
\end{itemize}

It is easily seen that $( \mathbf n, \mathbf r, \boldsymbol\psi)$ is a sufficient statistic for $\boldsymbol\theta$ based on the full data, just as $( \mathbf n, \mathbf r)$ is a sufficient statistic for $\boldsymbol\theta$ when we are only given the actually observed data.

Because of the constraints ${\bf C1},{\bf C2}$, we must have $\sum_{j=2}^J r_j\le K$, i.e. the number of non-singleton species observed is not larger than $K$.

Recall that $T=\sum_{j=1}^J r_j n_j$ is the total size of the sample. For given $\boldsymbol\psi$, define 
\begin{equation} \label{MH1}
n_0=r_1-\sum_{a=1}^K1\{\boldsymbol\psi(a)=1\},  %\eqno(1)
\end{equation}
the total number of times a blob species was observed. The full data likelihood is
\begin{equation}
\frac{T!}{n_0! \prod_{1\le a\le K : \boldsymbol\psi(a)\ge 1} n_{\boldsymbol\psi(a)}!} ~~~ p_0^{n_0}  \prod_{1\le a \le K\,:\,\boldsymbol\psi(a)\ge 1} p^{n_{\boldsymbol\psi(a)}}_a, \label{MH2} %\eqno(2)
\end{equation}
which, since the product over $a$ in the denominator of the multinomial term is equal to $\prod_{j=1}^J (n_j!)^{r_j}$, which is a constant, is proportional to (as a function of $\psi$ and $\theta$) 
\begin{equation}
\frac{1}{n_0!}~~~ p_0^{n_0}  \prod_{1\le a \le K\,:\,\boldsymbol\psi(a)\ge 1} p^{n_{\boldsymbol\psi(a)}}_a.\label{MH3} %\eqno(3)
\end{equation}

The observed data likelihood is the sum over all mappings $\boldsymbol\psi$ allowed by the constraints \textbf{C1} and \textbf{C2} of the full data likelihood. Note that $n_0$ occurs in the multinomial factor in the full data likelihood as well as as a power of $p_0$, and that $n_0$ depends on $\boldsymbol\psi$. 

{%In writing out (2) we have taken advantage of the facts that $0!=1$ and $p_a^0=1$ to omit $a$ such that $\boldsymbol\psi(a)=0$ from both of the products over non-blob population species. (KEEP IN TECH REPORT)}

\subsection{The moves}
We will define a random walk on the set of all mappings $\boldsymbol\psi$ allowed by the constraints \textbf{C1} and \textbf{C2}. It will be a Markov process with the set of mappings $\boldsymbol\psi$ as  the (huge) state-space; the graph of possible transitions between states will however be sparse. Inspection of the likelihood \eqref{MH3} suggests two kinds of moves: $(i)$ An \emph{exchange move}: exchanging the values of $\boldsymbol\psi(a)$ and $\boldsymbol\psi(b)$ for a chosen pair of different non-blob population species $a$ and $b$ such that $\boldsymbol\psi(a)\ne 0$, $\boldsymbol\psi(b)\ne 0$, and $\boldsymbol\psi(a)\ne\boldsymbol\psi(b)$, and $(ii)$ A \emph{blob move}: increasing or decreasing $n_0$ by one by choosing an $a$ such that $\boldsymbol\psi(a)=0$ or $\boldsymbol\psi(a)=1$ and exchanging the value $0$ of $\boldsymbol\psi(a)$ for $1$ or vice-versa.

Notice that these moves are not always possible. 
\blem
$(i)$: If $J\ge 3$ an exchange move is always possible. $(ii)$: If $S>0$, where $S=r_1$ is the number of singletons, and $K>N$, where $N=\sum_{j=2}^J r_j$ is the number of non-singletons, a blob move is always possible. 
\elem
\prf$(i)$: If we cannot find distinct $a$, $b$ with $\boldsymbol\psi(a)\ne 0$ , $\boldsymbol\psi(b)\ne 0$, and  $\boldsymbol\psi(a)\ne \boldsymbol\psi(b)$, an exchange move is impossible. However, as long as $J\ge 3$ there are at least two non-blob species observed a different number of times, and an exchange move is always possible. 

$(ii)$: It is always possible \emph{either} to increase or to decrease $n_0$ but it is not always possible to do both, since there is a minimum value, which can only be increased, and a maximum value, which can only be decreased (unless the minimum and maximum possible values of $n_0$ coincide). 

The maximum possible value of $n_0$, the number of times a blob species is observed, is the number of singletons  $S$ in the sample, and it is feasible to let every singleton correspond to a blob species. 

To determine the minimal value, define $L=N+S$, the length of the observed partition of $T$. The number of population species $a$ associated by $\boldsymbol\psi$ with singletons, i.e. such that  $\boldsymbol\psi(a)=1$, cannot exceed the total number of singletons $S$ but it also cannot exceed $K-N$. It can equal the minimum of these two numbers. Thus the lower bound on $n_0$ is given by the requirement $S-n_0\le \min(S,K-N)$, which is equivalent to $-n_0 \le \min(0,K-L)$, which is equivalent to $n_0\ge \max(0,L-K)$.

In summary, $\max(0,L-K)\leq n_0\leq S$ and therefore as long as $\max(0,L-K)<S$ or equivalently $S>0$ and $L-K<S$, thus $K>N$, a blob move is always possible. 
\eop

An exchange move defined by choice of a pair $(a,b)$ is its own reverse; and a blob move defined by choice of a single $a$ is its own reverse too. Moreover the number of candidate pairs $(a,b)$ for an exchange move is the same before and after the move. The number of candidates $a$ for a blob move is also the same before and after the move, except perhaps when $n_0$ is minimal or maximal. We shall further investigate these extreme cases later.

\subsubsection{Exchange moves}

For an exchange move we pick uniformly at random distinct $a$ and $b$ such that $\boldsymbol\psi(a)\ne 1$, $\boldsymbol\psi(b)\ne 1$, $\boldsymbol\psi(a)\ne\boldsymbol\psi(b)$.
The Metropolis factor follows from the formula \eqref{MH3} for the full data likelihood. The move would convert the factor $p_a^{n_{\boldsymbol\psi(a)}}p_b^{n_{\boldsymbol\psi(b)}}$ into  $p_a^{n_{\boldsymbol\psi(b)}}p_b^{n_{\boldsymbol\psi(a)}}$. The logarithm of the ratio of the full data likelihood ``after'' to ``before'' equals 
\begin{eqnarray*}
&(n_{\boldsymbol\psi(a)}\log p_b + n_{\boldsymbol\psi(b)}\log p_a) - (n_{\boldsymbol\psi(a)}\log p_a + n_{\boldsymbol\psi(b)} \log p_b) &\\
&= (n_{\boldsymbol\psi(a)} -n_{\boldsymbol\psi(b)})(\log p_b - \log p_a).&
\end{eqnarray*}
Thus we draw $Z$ from the standard exponential distribution and accept the move if and only if, since $\exp(-Z)$ is $\mathrm{Unif}[0,1]$-distributed,
\begin{equation}
- Z ~ \le ~ (n_{\boldsymbol\psi(a)} -n_{\boldsymbol\psi(b)})(\log p_b - \log p_a).\label{MH4}%\eqno(4)
\end{equation}
If the right hand side of \eqref{MH4} is positive, its exponent is larger than 1, and the move is accepted. If  the right hand side of \eqref{MH4} is negative, its exponent lies between $0$ and $1$, and hence the move is accepted with probability equal to this exponent.

\subsubsection{Blob moves}
In order to describe a blob move we separate between the three cases where (i) $n_0$ is equal to its minimal value, $\max(0,L-K)$, or (ii) maximal value, $S$, or (iii) is somewhere in between. 

If $n_0=\max(0,L-K)$, we pick a population species uniformly at random from the set $\{a:\boldsymbol\psi(a)=1\}$. If $n_0=S$, we pick a population species uniformly at random from the set $\{a:\boldsymbol\psi(a)=0\}$. When neither extreme case holds, we pick a sample species uniformly at random from the set $A(\boldsymbol\psi)=\{a:\boldsymbol\psi(a)=0\mbox{ or }\boldsymbol\psi(a)=1\}$.

However when $n_0=S$, there actually are no $a$ with $\boldsymbol\psi(a)=1$, so the rule prohibiting us to pick one of such $a$ in this case is superfluous. Similarly, if $n_0=L-K \ge 0$ then there are no $a$ with  $\boldsymbol\psi(a)=0$, and again the prohibition on picking such $a$ in this case is superfluous. Thus the rule for picking $a$ is simpler than first appeared: We always pick a population species uniformly at random from the set $A(\boldsymbol\psi)$. The number of species in $A(\boldsymbol\psi)$ is $K-N$, except when $L<K$ and $n_0=0$, in which case $A(\boldsymbol\psi)=\{a:\boldsymbol\psi(a)=1\}$ and then the number of species is $S < K-N$.

After the random choice of a species $a$ from $A(\boldsymbol\psi)$, the proposed move is to exchange the value of $\boldsymbol\psi(a)$ from $0$ to $1$ or vice-versa. In the first case $n_0$ is decreased by one, a factor $p_a^1$ gets added to the product of probabilities in \eqref{MH3}, and the logarithm of the Metropolis contribution to the acceptance criterion is $\log p_a -\log p_0  +\log n_0$. In the second case $n_0$ is increased by one, a factor $p_a^1$ gets deleted, and the logarithm of the Metropolis contribution to the acceptance criterion is  $\log p_0 -\log p_a -\log (n_0+1)$.

Recall that the Hastings factor in the Metropolis-Hastings algorithm is the ratio of the probabilities of the reverse move to the forward move. We have seen that, with one exception, the number of choices for $a$ is equal, both before and after the move, to $K-N$, so in general there is no Hastings contribution. The exceptional case is when $L<K$, and $n_0=0$ and is about to be increased by $1$ (because we picked $a$ with $\psi(a)=1$), or $n_0=1$ and is about to be decreased by $1$ (because we picked $a$ with $\psi(a)=0$). In these two cases the number of choices for the forward move and the reverse move are $S$ and $K-N>S$, and vice versa.

This means that if for a blob move we have picked $a$ with $\boldsymbol\psi(a)=1$, the move (put $\boldsymbol\psi(a)=0$) is accepted if and only if
\begin{equation}%$$ 
-Z ~ \le ~ \log p_0 - \log p_a -\log (n_0+1). %\eqno(5)$$
\end{equation}
(where $Z$ is a standard exponential random variable),
\emph{except} when $L < K$ and $n_0=0$, when the acceptance criterium is
\begin{equation}%$$ 
-Z ~ \le ~  \log p_0 - \log p_a - \log(K-N) + \log S %\eqno(6)$$
\end{equation}

If on the other hand we have picked $a$ with $\boldsymbol\psi(a)=0$, the move (put $\boldsymbol\psi(a)=1$) is accepted if and only if
\begin{equation}%$$ 
-Z ~ \le ~ \log p_a - \log p_0 + \log n_0.%\eqno(7)$$
\end{equation}
\emph{except} when $L < K$ and $n_0=1$, when the acceptance criterium is
\begin{equation}%$$ 
-Z ~ \le ~  \log p_a - \log p_0 + \log(K-N) - \log S. %\eqno(8)$$
\end{equation}

\subsection{The  SA-EM}
We next describe the ``statistical part'' of the SA-EM algorithm. We use it to make a so called ``stochastic approximation'' of the conditional expectation of an underlying full data sufficient statistic given the actually observed data at the current parameter estimates, in the E step, and then to re-estimate the parameters by maximum likelihood using the current approximation of the full data sufficient statistic, in the M step.

A suitable choice for the sufficient statistic is the vector $\mathbf g$ of relative frequencies $g_a$, $0\le a\le K$, of the underlying population species in our sample of size $T$. Given the vector $\mathbf f$ of observed sample species distinct relative frequencies $f_j=n_j/T$, $1 \le j \le J$, and given a realisation of the ``missing'' map $\boldsymbol\psi$, the underlying population relative frequencies are uniquely determined, for $1\le a \le K$,
\begin{eqnarray*}
g_a &=& f_j,\quad\text{ if }\boldsymbol\psi(a)=j\ge 1,\\
\quad g_a &=& 0, \quad\text{ if }\boldsymbol\psi(a)=0,\\
g_0 &=& \frac{n_0}{T}.
\end{eqnarray*}

At any point in the iterations we have a running estimate, denoted by $\boldsymbol\mu = (\mu_a)_{0\le a\le K}$, of the \emph{conditional expectation} of the vector $\mathbf g$ given the observed data $(\mathbf n,\mathbf r)$. The expectation is taken under the current estimate of the vector of probabilities $\boldsymbol\theta=(p_a)_{0\le a\le K}$. We generate a new \emph{realisation} of $\mathbf g$ as just defined, thought to be a realisation from the distribution of $\mathbf g$ given $(\mathbf n,\mathbf r)$ again under the current estimate of the parameter $\boldsymbol\theta$. 

The stochastic approximation update is to replace the current estimate $\mathbf m$ of the conditional expectation of $\mathbf g$ given $(\mathbf n,\mathbf r)$ under $\boldsymbol\theta$ by a weighted average of its current value and the current realisation $\mathbf g$ drawn from the distribution of $\mathbf g$ given $(\mathbf n,\mathbf r)$ under $\boldsymbol\theta$:  replace $\boldsymbol\mu$ by $(1-\gamma)\boldsymbol\mu + \gamma \mathbf g$ where the weight $\gamma$ will be a function of the iteration number, which we  denote by $k$. 

These weights should satisfy $\sum_k \gamma_k=\infty$, $\sum_k \gamma_k^2 < \infty$, cf. \cite{Delyon99}. Many authors propose to take $\gamma_k=1/k$ but we found that $\gamma_k=1/k^{2/3}$ worked better. For small test problems, we found that an initial value of $k=k_0=1000$ gave good results in conjunction with $\gamma_k=1/k^{2/3}$. 

In the E-step we take the current value of the vector $\boldsymbol\mu=(\mu_a)_{0\le a\le K} $ and maximise the log likelihood $\sum_{0\le a\le K} \mu_a \log p_a$ subject to the constraints $p_1\ge p_2\ge \dots p_K\ge 0$, $p_0\ge 0$, $\sum_{0\le a\le K} p_a=1$. This is equivalent to taking $p_1\ge \dots \ge p_K$ as a (version) of the isotonic (decreasing) regression of the vector $(\mu_a: 1\le a\le K)$ which can be found using a modification of the well-known \emph{pool adjacent violators} algorithm, supplemented with the assignment $p_0=\mu_0$, see \cite{robertson:wright:dykstra:1988} and the comments in the discussion. % see the Appendix for details.

Apart from the initialisation of $k$, also a realization of the mapping $\boldsymbol\psi$, an value of $\boldsymbol\theta$, and a  value of $\boldsymbol\mu$ need to be initialized. Since at convergence of the algorithm, $\boldsymbol\mu$ and $\boldsymbol\theta$ will be equal to one another, it also makes sense to initialise them equal to one another. A neutral initial guess for $\theta$ would be a  defective uniform probability distribution on $\{1,\dots,K\}$ supplemented with a not too small positive mass $p_0$ for the blob.  

We initialise $\boldsymbol\psi$, thought of as a vector, by assigning its first $r_J$ components all with the value $J$, the next $r_J-1$ with the value $J-1$, and so on, until we get to the $r_2$ components assigned with the value $2$; all remaining components are assigned the value $0$. Thus, under $\boldsymbol\psi$, a more frequently observed sample species has a larger population probability than a less frequently observed sample species, and all singletons are actually blob species.

% ++++++++++++++++++++++++++++++++++++++++++++++++++++++++++++++++++++++++++
% ++++++++++++++++++++++++++++++++++++++++++++++++++++++++++++++++++++++++++
% ++++++++++++++++++++++++++++++++++++++++++++++++++++++++++++++++++++++++++
% ++++++++++++++++++++++++++++++++++++++++++++++++++++++++++++++++++++++++++
% ++++++++++++++++++++++++++++++++++++++++++++++++++++++++++++++++++++++++++

% \input section_lowerbound
%
 \section{An algorithm for estimating a decreasing multinomial probability with lower bound} \label{App:3}
 In this appendix we present an algorithm for bounded isotonic regression and prove its convergence.
 
 \subsection{The algorithm and its convergence}
 
 Assume we have observations $(x_1,\ldots,x_n)$ of a multinomial random variable $\mathrm{Multi}(n,p)$ where $n=\sum_{i=1}^k x_i$ with $p=(p_1,\ldots,p_k)$ a vector a numbers $p_i\in [0,1]$ such that $\sum_{i=1}^k p_i=1$ and $p_1\geq \ldots \geq p_k$. Assume that the vector $p$ is unknown and assume also that for a given constant $0<c<1/k$ we know that $p_k\geq c$. The goal is then to estimate $p$ under the assumption that 
\begin{eqnarray}\label{eq:p_restriction_1}
A^{(c)}(p)&=&\{p_1\geq \ldots \geq p_k\geq c\}
\end{eqnarray}   
holds. 
Note that the restrictions $(\ref{eq:p_restriction_1})$ can be written as
\begin{eqnarray} \label{eq:p_restriction_2}
 A^{(c)}(p)&=&\cup_{j=1}^k A_j^{(c)}(p)
\end{eqnarray}
with
\begin{eqnarray}  
  A_j^{(c)}(p)&=&\{p_1\geq \ldots \geq p_j\} \cap \{p_{j+1}=\ldots=p_k\}\cap \{p_j\geq c\}, \label{eq:p_restriction_3}\\
A_{j}^{(c)}(p)&\subset &A_{j+1}^{(c)}(p),\mbox{ for all $p,j$ and }c.\label{eq:p_restriction_4}
\end{eqnarray}

Let ${\cal F}_c=\{q\in [0,1]^k: \sum_{i=1}^k q_i=1, q_1\geq \ldots \geq q_k\geq c\}$. Define the likelihood and log likelihood as
\begin{eqnarray*}
   L(p)&=&\frac{n!}{x_1! \cdots x_k!} p_1^{x_1}\cdots p_k^{x_k}\\
   \log L(p)&\sim&\sum_{i=1}^n x_i\log(p_i)=:l(p)
\end{eqnarray*}
 and define the order restricted mle $\hat{p}^{(c)}$ with lower bound at $c$ as (where appropriate we suppress the explicit dependence of the estimate on $c$ in the notation, and thus write $ \hat{p}= \hat{p}^{(c)}$)
 \begin{eqnarray*}
      \hat{p}^{(c)}=(\hat{p}_1,\ldots,\hat{p}_k)&=&\argmax_{q \in {\cal F}_c} \sum_{i=1}^k x_i \log q_i =\argmax_{q \in {\cal F}_c}  l(q).
 \end{eqnarray*}

The linear restriction $\sum_{i=1}^k q_i=1$ can be taken care of by introducing a Lagrange multiplier, so that the optimization is equivalent to maximization of
\begin{eqnarray*}
     \tilde{l}(q,\lambda)&=& \sum_{i=1}^n x_i\log(q_i)- \lambda (\sum_{i=1}^{k} q_i -1),
\end{eqnarray*}
with respect to $q$, over the set ${\cal F}_c  =\{q\in [0,1]^k:  q\in A^{(c)}(q)\}$. 

Note that with ${\cal F}_{c,j}  = \{q\in [0,1]^k:  q\in A_j^{(c)}(q)\}$, by $(\ref{eq:p_restriction_1}),(\ref{eq:p_restriction_2}),(\ref{eq:p_restriction_3})$ we obtain
\begin{eqnarray}
   {\cal F}_c  &=&\cup_{j=1}^{k}{\cal F}_{c,j}  ,\label{eq:set_restriction_1} \\
  {\cal F}_{c,j}   &\subset& {\cal F}_{c,j+1}  , \mbox{ for all }c,j.\label{eq:set_restriction_2} 
\end{eqnarray}
This  shows that (since ${\cal F}_{c,k}  ={\cal F}_{c}  $, this is  only the definition of the MLE again) 
\begin{eqnarray*}
    \hat{p}^{c}&=&\argmax_{q \in  {\cal F}_{c,k}   , \lambda} \tilde{l}(q,\lambda).
\end{eqnarray*}

We will maximize $\tilde{l}$ over ${\cal F}_c  $, by going through the sets ${\cal F}_{c,k}  ,\ {\cal F}_{c,k-1}  \ldots$. Note that the function ${l}$ is strictly concave, since it is a composition of the strictly concave function $log$ and a linear function, and the set ${\cal F}_c$ over which one wants to optimize $l$ is a closed convex set, since ${\cal F}_c$ is an intersection of closed convex cones and a hyperplane. This implies that there is a global solution $\hat{p}^c$, that it is unique and that $\hat{p}^c$ lies in ${\cal F}_c$. 

We present an iterative algorithm for obtaining the solution; this goes via a sequence $ \hat{p}^{(1)}, \hat{p}^{(2)},\ldots, \hat{p}^{(j)} $, with $j \in \{1,\ldots,k\}$ a random index, and such that $\hat{p}^{(j)}=   \hat{p}^{c}$ so that the algorithm converges, in a finite number of steps. 

\begin{algo} (Bounded isotonic regression) \label{algo1}

\begin{enumerate}
\item  Let
\begin{eqnarray*}
          \hat{p}^{(1)}&=&\argmax_{q \in  {\cal F}_{0,k}   , \lambda} \tilde{l}(q,\lambda),
\end{eqnarray*}
be the unrestricted estimator. Then if $ \hat{p}^{(1)}_k\geq c$, we are done, and the algorithm has converged with $j=1$. 

\item If not, so if $ \hat{p}^{(1)}_k<c $, we define the next element $\hat{p}^{(2)}$ in the sequence as follows: Let
$\hat{p}^{(2)}_k= c$. Define
\begin{eqnarray*}
     (\hat{p}_1^{(2)},\ldots, \hat{p}_{k-1}^{(2)})&=&\argmax_{q_1\geq \ldots \geq q_{k-1}\geq 0, \lambda}  \tilde{l}(q_1,\ldots, q_{k-1},c,\lambda).
\end{eqnarray*}
If $\hat{p}_{k-1}^{(2)}\geq c$ we are done, and the algorithm has converged with $j=2$.

\item If not, so if $\hat{p}_{k-1}^{(2)}< c$, we define the next element $\hat{p}^{(3)}$ in the sequence by: Let $\hat{p}^{(3)}_{k-1}= \hat{p}^{(3)}_{k}=c$. Define
\begin{eqnarray*}
     (\hat{p}_1^{(3)},\ldots, \hat{p}_{k-2}^{(3)})&=&\argmax_{q_1\geq \ldots \geq q_{k-2}\geq 0, \lambda}  \tilde{l}(q_1,\ldots, q_{k-2},c,c,\lambda).
\end{eqnarray*}

\item Continued in the obvious way, until convergence. 
\end{enumerate}
\end{algo}

To prove that the algorithm converges we need the following simple result that has independent interest and is stated for completeness.

\blem\label{lem:concave}
Assume that $h$ is a function, concave over a convex set $I\subset {\mathbb R}^k$, with the set $I$ defined by inequalities and linear restrictions. Then, for any integer $0<b<k$,  the function
\begin{eqnarray*}
      g(x_{b+1},\ldots,x_k)&=&\sup_{(x_1,\ldots,x_b):(x_1,\ldots,x_b,x_{b+1},\ldots,x_k) \in I} h(x_1,\ldots,x_b,x_{b+1},\ldots,x_k)
\end{eqnarray*}
is concave over the convex set 
\begin{eqnarray*}
I_{\bar{b}}&=&\{(x_{b+1},\ldots,x_k)\in {\mathbb R}^b:(x_1,\ldots,x_b,x_{b+1},\ldots,x_k) \in I\}.
\end{eqnarray*}
\elem
\prf
By assumption $ I$ can be written as an intersection of sets of the form $\{x=(x_1,\ldots,x_k): x_j\leq x_l\}$ (closed convex cones) and  $\{x=(x_1,\ldots,x_k):\sum_{i=1}^k a_i x_i=c\}$ (a hyperplane) for some real numbers $a_i, i=1,\ldots,k$ and $c$. Recall that by definition  $h$ is concave over $I$ if it's epigraph is a convex set. We use this repeatedly: 

The individual maximization $h$ as above means projecting the epigraph of $h$ on the faces $\{x=(x_1,\ldots,x_k): x_j\leq x_l\}$. Each such projection (i.e. maximization) gives a convex set (in ${\mathbb R}^k$) and the intersection of them is convex, and thus the epigraph of a concave function over $I_{\bar{b}}$, which shows that $g$ is concave. \eop

\bth
If $0<c<1/k$, the algorithm above converges to the global maximum $\hat{p}^{(c)}$.
\eth
\prf
We want to maximize $\tilde{l}$ over ${\cal F}_c  $, which we do  by stepwise going through the sets ${\cal F}_{c,1}  ,  {\cal F}_{c,2}  ,\ldots, {\cal F}_{c,k}  $ starting with ${\cal F}_{c,k}  $. 

If after step $(i)$, $\hat{p}^{(1)}_k\geq c$, then since $ {\cal F}_{c}  \subset {\cal F}_{0}  $ the maximum over ${\cal F}_{0}  $ (i.e. $\hat{p}^{(1)}_k$) is equal to the maximum over ${\cal F}_{c}  $ and we are done. 

Assume that we are not done and instead that $\hat{p}^{(1)}_k< c$. Define the function 
\begin{eqnarray*}
    g_1(u)&=&\max_{q_1\geq \ldots \geq q_{k-1}\geq u, \lambda} \tilde{l}(q_1,\ldots,q_{k-1},u,\lambda).
\end{eqnarray*}
By the above lemma  $g_1$ is a concave function. Thus it attains it's maximum at some point, which is $\hat{p}_k^{(1)}$. From the assumption $c> \hat{p}^{(1)}_k$ and by the concavity of $g_1$ we have that $g_1(u)< g_1(c)$ for any $u>c> \hat{p}^{(1)}_k$, i.e. that for
\begin{eqnarray*}
     \max_{q_1\geq \ldots \geq q_{k-1}\geq u, \lambda} \tilde{l}(q_1,\ldots,q_{k-1},u,\lambda)  &<& \max_{q_1\geq \ldots \geq q_{k-1}\geq c, \lambda} \tilde{l}(q_1,\ldots,q_{k-1},c,\lambda).
\end{eqnarray*}
This implies that,  when maximizing under the restriction ${\cal F}_{c,k}  $, we can actually fix $u=c$ and look for
\begin{eqnarray*}
\argmax_{q_1\geq \ldots \geq q_{k-1}\geq c, \lambda} \tilde{l}(q_1,\ldots,q_{k-1},c,\lambda).
\end{eqnarray*}

We do this as follows: In step $(ii)$ we replace the lower bound $ q_{k-1}\geq c$ by $q_{k-1}\geq 0$. Thus we are looking over a larger set. If $\hat{q}^{(2)}_{k-1}\geq c$ we are done. If not, so if $c>\hat{q}^{(2)}_{k-1}$, we study the function
\begin{eqnarray*}
    g_2(u)&=&\max_{q_1\geq \ldots \geq q_{k-1}\geq 0, \lambda} \tilde{l}(q_1,\ldots,q_{k-1},c,\lambda).
\end{eqnarray*} 
An application of Lemma \ref{lem:concave} shows that $g_2$ is a concave function, it's maximum is attained at $\hat{q}^{(2)}_{k-1}$, and by concavity $g_2(u)< g_2(c)$ for any $u<c<\hat{q}_{k-1}^{(2)}$. Thus when maximizing under the restriction ${\cal F}_{c,k-1}  $, we may fix $u=c$ and look for
\begin{eqnarray*}
\argmax_{q_1\geq \ldots \geq q_{k-2}\geq c, \lambda} \tilde{l}(q_1,\ldots,q_{k-2},c,c,\lambda).
\end{eqnarray*}

We do this as follows: In step $(iii)$, replace the lower bound $ q_{k-2}\geq c$ by $q_{k-2}\geq 0$, and so on.

$(iv)$ Continue until convergence.

This scheme is clearly finite, and will terminate with $\hat{q}_{k-j}^{(j+1)}\geq c$ for some $j<k$, since $c<1/k$, $p_1\geq \ldots \geq p_k$ and $\sum_{i=1}^k p_k=1$.

Furthermore, the scheme searches for maxima, over the sets ${\cal F}_{c,k}   \setminus {\cal F}_{c,k-1},$ ${\cal F}_{c,k-1}   \setminus {\cal F}_{c,k-2}  , \ldots$. If it stops at an index $j$, that means that $\hat{p}_{k-j}^{(j+1)}\geq c$ and that $\hat{p}_{k-j}^{(j+1)}$ is obtained as a maximum over the set ${\cal F}_{0,k-j}   $, which is also a maximum over the set ${\cal F}_{c,k-j}   $, and we will have searched over the set 
\begin{eqnarray*}
\left (\cup_{i=1}^{j}  {\cal F}_{c,k-i+1}   \setminus {\cal F}_{c,k-i}   \right )\cup {\cal F}_{c,k-j}  &=&{\cal F}_c  ,
\end{eqnarray*}
i.e. the solution is a global solution. This ends the proof of convergence. \eop

  \subsection{Implementation of the algorithm}
  
  We now present a numerical implementation of Algorithm \ref{algo1} in C/C++. Recall that in the bounded isotonic regression $p_1\geq \ldots \geq p_k\geq c$. The algorithm uses the standard isotonic regression with  $p_1\geq \ldots \geq p_k\geq 0$, whose implementation is called \verb"isoreg( )" and can be found in the literature \cite{robertson:wright:dykstra:1988} (see also its implementation in {\tt R}). The implementation of Algorithm \ref{algo1} thus reads:
  
  \begin{lstlisting}
vector<double> isobound(vector<double> y, double c) {
    vector<double> isoreg(vector<double> y);
    int i = y.size();
    int j;
    double S = 0.0;
    for (j = 0; j < i; j++) {
        S += y.at(j);
    	} 
    double T = S; 
    vector<double> yf = isoreg(y);
    vector<double> z = y;
    while ( yf[i-1] < c ) {
        S = S-y[i-1];
        T = T-c;
        yf[i-1] = c;
        i = i-1;
        z.resize(i);
        vector <double> zf =isoreg(z);
        for (j=0;j<i;j++) yf[j] = zf[j]*T/S;
    	}
    return(yf);
}
  \end{lstlisting}

%+++++++++++++++++++++++++++++++++++++++++++++++++++++++++++++++++++++++++
%   References
%+++++++++++++++++++++++++++++++++++++++++++++++++++++++++++++++++++++++++
%\section*{References}
\addcontentsline{toc}{section}{References}

%\bibliography{hiprofile}

\begin{thebibliography}{10}

\bibitem{acharya:orlitsky:pan:2009}
J.~Acharya, A.~Orlitsky, and S.~Pan.
\newblock The maximum likelihood probability of unique-singleton, ternary, and
  length-7 patterns.
\newblock In {\em IIEEE International Symposium on Information Theory}, pages
  1135 -- 1139, 2009.

\bibitem{anevski:fougeres:2007}
D.~Anevski and A-L. Foug{\`e}res.
\newblock Limit properties of the monotone rearrangement for density and
  regression function estimation.
\newblock arxiv:0710.4617v1, Lund University, 2007.

%\bibitem{anevski:gill:zohren:2013}
%D.~Anevski, R.D. Gill and S. Zohren.
%\newblock Supplement to ``''Estimating a probability mass function with unknown labels``''.
%\newblock arxiv:1312.1200v2, 2013.


\bibitem{balabdaoui:rufibach:santambrogio:2009}
F. Balabdaoui, K. Rufibach, and F. Santambrogio.
\newblock Least squares estimation of two ordered monotone regression curves.
\newblock {\em Journal of Nonparametric Statistics}, 22:1019, 2009.

\bibitem{Delyon99}
B.~Delyon, M.~Lavielle, and E.~Moulines.
\newblock Convergence of a stochastic approximation version of the {EM}
  algorithm.
\newblock {\em Ann. Statist.}, 27(1):94--128, 1999.

\bibitem{dvoretzky:kiefer:wolfowitz:1956}
A.~Dvoretzky, J.~Kiefer, and J.~Wolfowitz.
\newblock Asymptotic minimax character of the sample distribution function and
  of the classical multinomial estimator.
\newblock {\em Ann. Math. Statist.}, 27:642--669, 1956.

\bibitem{efron:thisted:1976}
B.~Efron and R.~Thisted.
\newblock Estimating the number of unseen species: How many words did
  shakespeare know?
\newblock {\em Biometrika}, 63:435--447, 1976.

\bibitem{esty}
W.~W. Esty.
\newblock Confidence intervals for the coverage of low coverage samples.
\newblock {\em Ann. Statist.}, 10:190, 1982.

\bibitem{esty2}
W.~W. Esty.
\newblock A normal limit law for a nonparametric estimator of the coverage of a
  random sample.
\newblock {\em Ann. Statist.}, 11:905, 1983.

\bibitem{fisher:corbet:williams:1943}
R.A. Fisher, A.S. Corbet, and C.B. Williams.
\newblock The relation between the number of species and the number of
  individuals in a random sample of an animal population.
\newblock {\em J. Anim. Ecol.}, 12, 1943.

\bibitem{good:1953}
I.J. Good.
\newblock The population frequencies of species and the estimation of
  population parameters.
\newblock {\em Biometrika}, 40:237--264, 1953.

\bibitem{good:toulmin:1956}
I.J. Good and G.H. Toulmin.
\newblock The population frequencies of species and the estimation of
  population parameters.
\newblock {\em Biometrika}, 43:45--63, 1956.

\bibitem{hardy:littlewood:polya:1952}
G.~H. Hardy, J.~E. Littlewood, and G.~P{\'o}lya.
\newblock {\em Inequalities}.
\newblock Cambridge, at the University Press, 1952.
\newblock 2d ed.

\bibitem{jankowski:wellner:2009}
H.~Jankowski and J.A. Wellner.
\newblock Estimation of a discrete monotone distribution.
\newblock {\em Electron J Stat.}, 3:1567--1605, 2009.

\bibitem{lieb:loss:1996}
E.~H. Lieb and M.~Loss.
\newblock {\em Analysis}, volume~14 of {\em Graduate Studies in Mathematics}.
\newblock American Mathematical Society, 1996.

\bibitem{mao}
C.~X. Mao and B.~G Lindsay.
\newblock A {P}oisson model for the coverage problem with a genomic
  application.
\newblock {\em Biometrika}, 89:669, 2002.

\bibitem{massart:1990}
P.~Massart.
\newblock The tight constant in the {D}voretzky-{K}iefer-{W}olfowitz
  inequality.
\newblock {\em Ann. Probab.}, 18(3):1269--1283, 1990.

\bibitem{pan2009}
A.~Orlitsky and S.~Pan.
\newblock The maximum likelihood probability of skewed patterns.
\newblock In {\em IEEE International Symposium on Information Theory}, 2009.

\bibitem{orlitsky:sajama:santhanam:viswanathan:2004b}
A.~Orlitsky, S.~Sajama, N.P. Santhanam, K.~Viswanathan, and J.~Zhang.
\newblock Algorithms for modeling distributions over large alphabets.
\newblock In {\em Information Theory, 2004. ISIT 2004. Proceedings.
  International Symposium on Information Theory}, page 304, 2004.

\bibitem{orlitsky:sajama:santhanam:viswanathan:2004a}
A.~Orlitsky, S.~Sajama, N.P. Santhanam, K.~Viswanathan, and Junan Zhang.
\newblock On modeling profiles instead of values.
\newblock In {\em Proceeding UAI '04 Proceedings of the 20th conference on
  Uncertainty in artificial intelligence}, pages 426--435, 2004.

\bibitem{orlitsky:sajama:santhanam:viswanathan:2005}
A.~Orlitsky, S.~Sajama, N.P. Santhanam, K.~Viswanathan, and Junan Zhang.
\newblock Convergence of profile based estimators.
\newblock In {\em Information Theory, 2005. ISIT 2005. Proceedings.
  International Symposium on Information Theory}, pages 1843--1847, 2005.

\bibitem{orlitsky:santhanam:viswanathan:zhang:2004}
A.~Orlitsky, N.P. Santhanam, K.~Viswanathan, and Junan Zhang.
\newblock On modeling profiles instead of values.
\newblock In {\em Proceedings of the Twentieth Conference Annual Conference on
  Uncertainty in Artificial Intelligence (UAI-04)}, pages 426--435, Arlington,
  Virginia, 2004. AUAI Press.

\bibitem{ramanujan:hardy:1918}
S.~Ramanujan and G.H. Hardy.
\newblock Asymptotic formulae in combinatorial analysis.
\newblock {\em Proc. London Math. Soc.}, 17(1):75--115, 1918.

\bibitem{robertson:wright:dykstra:1988}
T.~Robertson, F.T. Wright, and R.L. Dykstra.
\newblock {\em Order Restricted Statistical Inference}.
\newblock John Wiley \& Sons Inc., New York, 1988.

\bibitem{vaneeden:1957}
C.~van Eeden.
\newblock Maximum likelihood estimation of partially or completely ordered
  parameters. ii.
\newblock {\em Proceedings Koninklijke Nederlandse Akademic van Wetenschappen,
  Series A. 60. Indagationes Mathematical}, 19:201--211, 1957.

\bibitem{vontobel:2012}
P.O.~Vontobel,
\newblock The Bethe Permanent of a Non-Negative Matrix.
\newblock {\em Information Theory, IEEE Transactions on  (Volume:59 ,  Issue: 3 )}, pages 1866 - 1901, 2012

\bibitem{vontobel:2014}
P. O.~Vontobel,
\newblock The Bethe and Sinkhorn approximations of the pattern maximum likelihood estimate and their connections to the Valiant-Valiant estimate.
\newblock {\em Proceedings of Information Theory and Applications Workshop (ITA)}, 9-14 Feb. 2014.


\bibitem{zhang:zhang:2009}
C-H. Zhang and Z.~Zhang.
\newblock Asymptotic normality of a nonparametric estimator of sample coverage.
\newblock {\em The Annals of Statistics}, 37:2582--2595, 2009.

\end{thebibliography}

\begin{thebibliography}{10}


\bibitem{Delyon99}
B.~Delyon, M.~Lavielle, and E.~Moulines.
\newblock Convergence of a stochastic approximation version of the {EM}
  algorithm.
\newblock {\em Ann. Statist.}, 27(1):94--128, 1999.


\bibitem{pan2009}
A.~Orlitsky and S.~Pan.
\newblock The maximum likelihood probability of skewed patterns.
\newblock In {\em IEEE International Symposium on Information Theory}, 2009.


\bibitem{robertson:wright:dykstra:1988}
T.~Robertson, F.T. Wright, and R.L. Dykstra.
\newblock {\em Order Restricted Statistical Inference}.
\newblock John Wiley \& Sons Inc., New York, 1988.


\end{thebibliography}

\end{document}